%% file: DDFG14.tex
\documentclass[11pt]{article}
\usepackage[dvips]{graphicx}
\usepackage{textcomp}
\usepackage{amsbsy}
\usepackage{latexsym}
\usepackage[mathscr]{eucal}
\usepackage{amsfonts}
\usepackage{amssymb}
\usepackage[usenames]{color}
\usepackage{fullpage}
\def\nnew{}
\def\mnew{}
\begin{document}

\title {Iteratively Re-weighted Least Squares Minimization for Sparse Recovery}

\author{Ingrid Daubechies\footnote{Princeton University, Department
of Mathematics and Program in Applied and Computational Mathematics,
{\tt ingrid@math.princeton.edu}.},~
Ronald DeVore\footnote{University of South Carolina, Department of
Mathematics, {\tt devore@math.sc.edu}.},~ 
Massimo Fornasier\footnote{Johann Radon Institute for Computational and 
Applied Mathematics, Austrian Academy of Sciences, 
{\tt massimo.fornasier@oeaw.ac.at}.},~ 
C.~Sinan G\" unt\" urk\footnote{New York University, 
Courant Institute of Mathematical Sciences, {\tt gunturk@courant.nyu.edu}.}}
\hbadness=10000
\vbadness=10000
\newtheorem{lemma}{Lemma}[section]
\newtheorem{prop}[lemma]{Proposition}
\newtheorem{cor}[lemma]{Corollary}
\newtheorem{theorem}[lemma]{Theorem}
\newtheorem{remark}[lemma]{Remark}
\newtheorem{example}[lemma]{Example}
\newtheorem{definition}[lemma]{Definition}
\newtheorem{proper}[lemma]{Properties}
\newtheorem{alg}{Algorithm}
%
\def\vp{\varphi}
\def\<{\langle}
\def\>{\rangle}
\def\t{\tilde}
\def\i{\infty}
\def\e{\varepsilon}
\def\sm{\setminus}
\def\nl{
newline}
\def\Chi{\raise .3ex \hbox{\large $\chi$}} 
\def\lsima{\hbox{\kern -.6em\raisebox{-1ex}{$~\stackrel{\textstyle<}{\sim}~$}}\kern -.4em}
\def\lsim{\hbox{\kern -.2em\raisebox{-1ex}{$~\stackrel{\textstyle<}{\sim}~$}}\kern -.2em}\def\[{\Bigl [}
\def\]{\Bigr ]}
\def\({\Bigl (}
\def\){\Bigr )}
\def\[{\Bigl [}
\def\]{\Bigr ]}
\def\({\Bigl (}
\def\){\Bigr )}
\def\L{�}
\def\pr{{\rm Prob}}
\newcommand{\iref}[1]{(\ref{#1})}
\def\ds{\displaystyle}
\def\ev#1{\vec{#1}}     
\newcommand{\lt}{\ell_{2}(\nabla)}
\def\Supp#1{{\rm supp\,}{#1}}
\def\e{\epsilon}
\def\R{\mathbb{R}}
\def\E{\mathbb{E}}
\def\nl{
newline}
\def\T{{\relax\ifmmode I\!\!\hspace{-1pt}T\else$I\!\!\hspace{-1pt}T$\fi}}
\def\N{\mathbb{N}}
\def\Z{\mathbb{Z}}
\def\Zd{\Z^d}
\def\Q{\mathbb{Q}}
\def\C{\mathbb{C}}
\def\Rd{\R^d}
\def\gsim{\mathrel{\raisebox{-4pt}{$\stackrel{\textstyle>}{\sim}$}}}
\def\sime{\raisebox{0ex}{$~\stackrel{\textstyle\sim}{=}~$}}
\def\lsim{\raisebox{-1ex}{$~\stackrel{\textstyle<}{\sim}~$}}
\def\div{\mbox{ div }}
\def\M{M}  \def\NN{N}                  
\def\L{{\ell}}               
\def\Le{{\ell_1}}            
\def\Lz{{\ell_2}}
\def\Let{{\tilde\ell_1}}     
\def\Lzt{{\tilde\ell_2}}
\def\Ltw{\ell_\tau^w(\nabla)}
\def\t#1{\tilde{#1}}
\def\la{\lambda}
\def\La{\Lambda}
\def\ga{\gamma}
\def\BV{{\rm BV}}
\def\Ga{\eta}
\def\al{\alpha}
\def\cZ{{\cal Z}}
\def\argmin{\mathop{\rm argmin}}
\def\argmax{\mathop{\rm argmax}}

\def\prob{\mathop{\rm prob}}

\def\cO{{\cal O}}
\def\cA{{\cal A}}
\def\cC{{\cal C}}
\def\cF{{\cal F}}
\def\bu{{\bf u}}
\def\bz{{\bf z}}
\def\bZ{{\bf Z}}
\def\bI{{\bf I}}
\def\cE{{\cal E}}
\def\cD{{\cal D}}
\def\cG{{\cal G}}
\def\cI{{\cal I}}
\def\cJ{{\cal J}}
\def\cM{{\cal M}}
\def\cN{{\cal N}}
\def\cT{{\cal T}}
\def\cU{{\cal U}}
\def\cV{{\cal V}}
\def\cW{{\cal W}}
\def\cL{{\cal L}}
\def\cB{{\cal B}}
\def\cG{{\cal G}}
\def\cK{{\cal K}}
\def\cS{{\cal S}}
\def\cP{{\cal P}}
\def\cQ{{\cal Q}}
\def\cR{{\cal R}}
\def\cU{{\cal U}}
\def\bL{{\bf L}}
\def\bK{{\bf K}}
\def\bC{{\bf C}}
\def\X{X\in\{L,R\}}
\def\ph{{\varphi}}
\def\D{{\Delta}}
\def\H{{\cal H}}
\def\bM{{\bf M}}
\def\bx{{\bf x}}
\def\bG{{\bf G}}
\def\bP{{\bf P}}
\def\bW{{\bf W}}
\def\bT{{\bf T}}
\def\bV{{\bf V}}
\def\bv{{\bf v}}
\def\bt{{\bf t}}
\def\bz{{\bf z}}
\def\bw{{\bf w}}
\def \span{{\rm span}}
\def \meas {{\rm meas}}
\def\rhom{{\rho^m}}
\def\lll{\langle}
\def\argmin{\mathop{\rm argmin}}
\def\argmax{\mathop{\rm argmax}}
\def\dJ{\nabla}
\newcommand{\ba}{{\bf a}}
\newcommand{\bb}{{\bf b}}
\newcommand{\bc}{{\bf c}}
\newcommand{\bd}{{\bf d}}
\newcommand{\bs}{{\bf s}}
\newcommand{\bff}{{\bf f}}
\newcommand{\bp}{{\bf p}}
\newcommand{\bg}{{\bf g}}
\newcommand{\by}{{\bf y}}
\newcommand{\br}{{\bf r}}
\newcommand{\be}{\begin{equation}}
\newcommand{\ee}{\end{equation}}
\newcommand{\bea}{$$ \begin{array}{lll}}
\newcommand{\eea}{\end{array} $$}
\def \Vol{\mathop{\rm  Vol}}
\def \mes{\mathop{\rm mes}}
\def \Prob{\mathop{\rm  Prob}}
\def \exp{\mathop{\rm    exp}}
\def \sign{\mathop{\rm   sign}}
\def \sp{\mathop{\rm   span}}
\def \vphi{{\varphi}}
\def \csp{\overline \mathop{\rm   span}}

%
\newcommand{\beqn}{\begin{equation}}
\newcommand{\eeqn}{\end{equation}}
\newcommand{\beqns}{$$}
\newcommand{\eeqns}{$$}
\def\beginproof{\noindent{\bf Proof:}~ }
\def\endproof{\hfill\rule{1.5mm}{1.5mm}\\[2mm]}

\newenvironment{Proof}{\noindent{\bf Proof:}\quad}{\endproof}

\renewcommand{\theequation}{\thesection.\arabic{equation}}
\renewcommand{\thefigure}{\thesection.\arabic{figure}}

\makeatletter
\@addtoreset{equation}{section}
\makeatother

\newcommand\abs[1]{\left|#1\right|}
\newcommand\clos{\mathop{\rm clos}\nolimits}
\newcommand\trunc{\mathop{\rm trunc}\nolimits}
\renewcommand\d{d}
\newcommand\dd{d}
\newcommand\diag{\mathop{\rm diag}}
\newcommand\dist{\mathop{\rm dist}}
\newcommand\diam{\mathop{\rm diam}}
\newcommand\cond{\mathop{\rm cond}\nolimits}
\newcommand\eref[1]{(\ref{#1})}
\newcommand\Hnorm[1]{\norm{#1}_{H^s([0,1])}}
\def\int{\intop\limits}
\renewcommand\labelenumi{(\roman{enumi})}
\newcommand\lnorm[1]{\norm{#1}_{\ell_2(\Z)}}
\newcommand\Lnorm[1]{\norm{#1}_{L_2([0,1])}}
\newcommand\LR{{L_2(\R)}}
\newcommand\LRnorm[1]{\norm{#1}_\LR}
\newcommand\Matrix[2]{\hphantom{#1}_#2#1}
\newcommand\norm[1]{\left\|#1\right\|}
\newcommand\ogauss[1]{\left\lceil#1\right\rceil}
\newcommand{\QED}{\hfill
\raisebox{-2pt}{\rule{5.6pt}{8pt}\rule{4pt}{0pt}}%
  \smallskip\par}
\newcommand\Rscalar[1]{\scalar{#1}_\R}
\newcommand\scalar[1]{\left(#1\right)}
\newcommand\Scalar[1]{\scalar{#1}_{[0,1]}}
\newcommand\Span{\mathop{\rm span}}
\newcommand\supp{\mathop{\rm supp}}
\newcommand\ugauss[1]{\left\lfloor#1\right\rfloor}
\newcommand\with{\, : \,}
\newcommand\Null{{\bf 0}}
\newcommand\bA{{\bf A}}
\newcommand\bB{{\bf B}}
\newcommand\bR{{\bf R}}
\newcommand\bD{{\bf D}}
\newcommand\bE{{\bf E}}
\newcommand\bF{{\bf F}}
\newcommand\bH{{\bf H}}
\newcommand\bU{{\bf U}}
\newcommand\cH{{\cal H}}
\newcommand\sinc{{\rm sinc}}
\def\enorm#1{| \! | \! | #1 | \! | \! |}

\newcommand{\dm}{\frac{d-1}{d}}

\let\bm\bf
\newcommand{\bbeta}{{\mbox{\boldmath$\beta$}}}
\newcommand{\bal}{{\mbox{\boldmath$\alpha$}}}
\newcommand{\bbi}{{\bm i}}

\newif\ifNZB
\newcommand\NZB[1]{\ifNZB \marginpar{\raggedright \scriptsize NZB:\\#1}
 \else \fi}
\newcommand{\FText}[1]{\mbox{#1}}
\makeatletter
\newcommand{\tr}{{\mathop{\operator@font T}\nolimits}}
\newcommand{\mod}{\mathop{\operator@font mod}}
\makeatother
\newcommand{\UArrow}[4]{
 \begin{array}{ll}
  #1&\stackrel{#2}\longrightarrow\\
  #3&\;\raisebox{1ex}{$\nearrow$}\mkern-14mu_{#4}
 \end{array}}
\newcommand{\DArrow}[4]{
 \begin{array}{ll}
  \stackrel{#1}\longrightarrow&#2\\
  \mkern-10mu_{#3}\mkern-22mu\raisebox{1ex}{$\searrow$}&#4
 \end{array}}
\newcommand{\fig}[3]{\par\begin{figure}[ht]
  \centerline{\epsfbox{#1.eps}}\caption{#3}\label{fig#2}\end{figure}}
\newcommand{\dI}{\Delta}
\maketitle
\date{}
\begin{abstract}
Under certain conditions (known as the {\em Restricted Isometry Property} or RIP) on the $m\times N$-matrix $\Phi$ (where $m<N$), vectors $x \in \R^N$ that are sparse (i.e. have most of their entries equal to zero) can be recovered exactly from $y:=\Phi x$ even though $\Phi^{-1}(y)$ is typically an $(N-m)$-dimensional hyperplane; in addition $x$ is then
equal to the element in $\Phi^{-1}(y)$ of minimal
$\ell_1$-norm. This minimal element can be identified via linear programming algorithms. \\
We study an alternative method of determining $x$, as the limit of an 
{\em Iteratively Re-weighted Least Squares} (IRLS) algorithm.
The main step of this IRLS finds, for a given weight vector $w$, the element
in $\Phi^{-1}(y)$ with smallest $\ell_2(w)$-norm.  If $x^{(n)}$ is the solution at iteration step $n$, then the new weight $w^{(n)}$ is defined by 
$w^{(n)}_i:=\left[|x^{(n)}_i|^2+\epsilon_n^2\right]^{-1/2}$, $i=1,\dots,N$, for a decreasing
sequence of adaptively defined $\epsilon_n$; this  updated weight is then used to obtain $x^{(n+1)}$ and the process is repeated. We prove that when $\Phi$
satisfies the RIP conditions, the sequence $x^{(n)}$ converges for all 
$y$, regardless of
whether $\Phi^{-1}(y)$ contains a sparse vector. If there is a sparse 
vector in $ \Phi^{-1}(y)$, then the limit is this sparse vector, and
when $x^{(n)}$ is sufficiently  close to the limit, the remaining steps of the algorithm  converge exponentially fast ({\em linear convergence} in the terminology of numerical optimization).  
The same algorithm with the ``heavier'' weight
$w^{(n)}_i = \left[|x^{(n)}_i|^2+\epsilon_n^2\right]^{-1+\tau/2}$, $i=1,\dots,N$, where $0<\tau<1$, can recover sparse solutions as well;
more importantly, we show its local convergence is superlinear and approaches a 
{\em quadratic} rate for $\tau$ approaching to zero.
\end{abstract}
\section{Introduction}
\label{intro}

Let $\Phi$ be an $m\times N$ matrix with $m<N$ and let $y\in \R^m$.
(In the compressed sensing application that motivated this study,
$\Phi$ typically has full rank, i.e. $\mbox{Ran}(\Phi)= \R^m$.
We shall implicitly assume, throughout the paper, that this is the case.
Our results still hold for the case where 
$\mbox{Ran}(\Phi)\subsetneq \R^m$, with the proviso that $y$ must
then lie in $\mbox{Ran}(\Phi)$.)

The linear system of equations 
\beqn
\label{system}
\Phi x=y
\eeqn
is underdetermined, and has infinitely many solutions.  If $\cN:=\cN(\Phi)$ is the null space of $\Phi$ and $x_0$ is any solution to \eref{system} 
then the set $\cF(y):=\Phi^{-1}(y)$ of all solutions to \eref{system} is given by $\cF(y)=x_0+\cN$.

In the absence of any other information, no solution to \eref{system} is to be preferred over any other.   However, many scientific applications work under the assumption that the desired solution $x\in\cF(y)$ is either sparse or well approximated by (a) sparse vector(s).  Here and later, we say a vector {\em has sparsity} $k$ (or is {\em $k$-sparse}) if it has {\it at most} $k$ nonzero coordinates. 
  Suppose then that we know that the desired solution of \eref{system} is $k$-sparse, where $k<m$ is known. 
How could we find such an $x$?   One possibility is to consider any set $T$
of $k$ column indices and find the least squares solution 
$x^T:=\argmin_{z\in \cF(y)}\|\Phi_Tz-y\|_{\ell_2^m}$,  where $\Phi_T$ is
obtained from $\Phi$ by setting to zero all entries that are not in columns from $T$.   Finding $x^T$ is numerically simple 
(see \eref{mfsol}).  After finding each $x^T$, we choose the particular 
set $T^*$ that minimizes the residual $\|\Phi_Tz-y\|_{\ell_2^m}$ .  This  
would find   a $k$-sparse solution (if it exists),  $x^*=x^{T^*}$. However, this 
naive method is numerically prohibitive when $N$ and $k$ are large, since it requires solving $  N \choose k$ least squares problems.

An attractive alternative to the naive minimization is 
its convex relaxation that consists in selecting the element in $\cF(y)$ which has minimal $\ell_1$-norm:
\beqn
\label{l1}
x:=\argmin_{z\in\cF(y)}\|z\|_{\ell_1^N}.
\eeqn
Here and later we use the $\ell_p$-norms
 \beqn
\label{lp}
\|x\|_{\ell_p}:=\|x\|_{\ell_p^N}:=\left\{\begin{array}{ll}
\left(\sum_{i=1}^N|x_j|^p \right)^{1/p}, &0 < p<\infty,\\
 \max_{j=1,\dots,N} |x_j|,& p=\infty.
\end{array}\right.
\eeqn

Under certain assumptions on $\Phi$ and $y$  that we shall describe in \S \ref{unsection},
it is known that \eref{l1} has a unique solution (which we shall denote by $x^*$), and that, when there is a $k$-sparse solution to \eref{system}, \eref{l1} will find this solution \cite{carota06-1,cata05,do04,do06-1}. 
Because the problem \eref{l1} can be formulated as a  linear program, 
it is numerically tractable.   

Solving underdetermined systems by $\ell_1$-minimization  has a long history.  It is at the heart of many numerical algorithms for approximation, compression, and statistical estimation.
The use of the $\ell_1$-norm as a sparsity-promoting functional can be found first in
reflection seismology and in deconvolution of seismic traces \cite{CM,SS,TBM}.
Rigorous results for $\ell_1$-minimization began to appear in the late-1980's, with Donoho and Stark \cite{dost89} and Donoho and Logan \cite{dolo92}. Applications for $\ell_1$-minimization in statistical estimation began in the mid-1990's with the introduction of the LASSO and related formulations \cite{ti96} (iterative soft-thresholding), also known  as
Basis Pursuit \cite{chdo98}, proposed in compression applications for extracting the sparsest signal representation from highly overcomplete frames. 
  Around the same time other signal processing groups started using $\ell_1$-minimization for the analysis of sparse signals; see, e.g. \cite{OBSK}.
The applications and understanding of 
$\ell_1$-minimization saw a dramatic increase in the last 5 years
\cite{do04,dota05,do06-1,dota06,cata05,cataXX,carota06-1, carotaXX}, with
 the development of fairly general mathematical 
frameworks in which
$\ell_1$-minimization, known heuristically to be sparsity-promoting, 
can be {\em proved} to recover sparse solutions {\em exactly}.
We shall not trace all the relevant results and applications;
a detailed history is beyond the scope of this introduction. 
We refer the reader to the survey papers \cite{ca06,ba07}. 
The reader can also find a  comprehensive collection of the ongoing recent developments at the web-site {\tt http://www.dsp.ece.rice.edu/cs/}.
In fact,  $\ell_1$-minimization has been so surprisingly effective in several applications, that Cand\`es, Wakin, and Boyd call it the ``modern least squares'' in \cite{cawaboXX}. We thus clearly need efficient algorithms for the minimization problem (\ref{l1}). 

Several alternatives to \eref{l1}, see, e.g.,  \cite{dots06-1,li93}, have been proposed as possibly more efficient numerically, or simpler to implement by non-experts, than standard algorithms for linear programming 
(such as interior point or barrier methods). 
In this paper we clarify fine convergence properties of one such 
alternative method, called  {\em Iteratively Re-weighted Least Squares minimization} (IRLS).  It begins with the following observation (see \S \ref{unsection} for details).  
If \eref{l1} has a 
solution $x^*$ that has no vanishing coordinates, then
the (unique!) solution $x^w$ of the weighted least squares problem
\beqn
\label{wls}
x^w:=\argmin_{z\in\cF(y)}\|z\|_{\ell_2^N(w)},\quad w:=(w_1,\ldots,w_N),
\quad \mbox{ where } w_j:=|x_j^*|^{-1},
\eeqn
coincides with $x^*$.
(The following argument provides a short proof by contradiction of this statement.
Assume that $x^*$ is not 
the $\ell_2^N(w)$-minimizer. Then there exists $\eta \in \cN$ such 
that $\| x^* + \eta \|_{\ell_2^N(w)}^2 < \| x^*\|_{\ell_2^N(w)}^2$ or equivalently $\frac{1}{2} \| \eta \|_{\ell_2^N(w)}^2 < - \sum_{j=1}^N w_j \eta_j x^*_j = \sum_{j=1}^N  \eta_j \sign(x^*_j) $. However, because $x^*$ is an $\ell_1$-minimizer, we have 
$\|x^*\|_{\ell_1} \leqslant  \|x^*+ h\eta\|_{\ell_1}$
for all $h\neq 0$; taking $h$ sufficiently small, this 
implies $\sum_{j=1}^N  \eta_j \sign(x^*_j) =0$, a contradiction.)

Since we do not know $x^*$, this observation cannot be used directly.
However, it leads to the following paradigm for finding $x^*$.   We choose a starting weight $w^0$ and solve \eref{wls} for this weight.
We then use this solution   to define a new weight $w^1$ and repeat this process.  
An IRLS algorithm of this type appears for the first time in the approximation practice in the Ph.D. thesis of Lawson in 1961 \cite{LW61}, in the form of an algorithm for solving uniform approximation problems, in particular by Chebyshev polynomials, by means of limits of weighted $\ell_p$--norm solutions. This iterative algorithm is now well-known in classical approximation theory as Lawson's algorithm. In \cite{cl72} it is proved that this algorithm has in principle a linear convergence rate. In the 1970s extensions of Lawson's algorithm for $\ell_p$-minimization, and in particular $\ell_1$-minimization, were proposed.  In signal analysis, IRLS was
proposed as a technique to build algorithms for sparse signal reconstruction
in \cite{GoRa}.  
Perhaps the most comprehensive   mathematical analysis of the performance of IRLS for 
$\ell_p$-minimization 
was given in the work of Osborne \cite{O}.

Osborne proves that a suitable IRLS method is convergent for $1 < p < 3$. 
For $p=1$, if $w^n$ denotes the weight at the $n$th iteration and $x^n$ 
the minimal weighted least squares solution for this weight, then the algorithm considered by Osborne defines the new weight $w^{n+1}$ coordinatewise as 
$w_j^{n+1}:=|x_j^n|^{-1}$.  His main conclusion in this case is that 
if the $\ell_1$ minimization problem \eref{l1} has a unique solution, then
the algorithm converges to this solution, in principle with linear convergence rate, i.e. exponentially fast, with a constant ``contraction
factor''.

However, the analysis of Osborne does not take into consideration what happens if one of the coordinates vanishes at some iteration $n$, i.e. $x_j^n=0$.
Taking this to impose that the corresponding weight component $w_j^{n+1}$ must ``equal'' $\infty$ leads to 
$x_j^{n+1}=0$ at the next iteration as well; this 
then persists in all later iterations. If $x^*_j=0$, all is well, but
if there is an index $j$ for which
$x^*_j\neq 0$, yet $x_j^n=0$ at some iteration step $n$, then this ``infinite weight'' prescription leads to problems. In practice, this is avoided by changing the definition of the weight at coordinates $j$ where $x_j^n=0$ (see \cite{li93} and \cite{ChL,FM} where a variant for total variation minimization is studied); such
modified algorithms need no longer converge to $x^*$, however).  
Because
Osborne's convergence proof is local, it implies that if the iterations
begin with a vector sufficiently close to the solution, and if the solution is unique and has only nonzero entries, then none of
the $x_j^n=0$ vanish, and the weight-change is not required; Osborne's analysis does indeed show the linear convergence rate of the algorithm
under these assumptions. Unfortunately, as we will see in Remark \ref{vanishing}, the uniqueness of the solution necessarily implies that it has vanishing components.  In other words, the set of vectors to which Osborne's analysis applies is vacuous.

The purpose of the present paper is to put forward an IRLS algorithm
that gives a re-weighting without infinite components in the weight, and to  provide an analysis of this algorithm, with various results about its convergence and rate of convergence. 
It turns out that care must be taken in just how the new weight $w^{n+1}$ is derived from the solution $x^{n}$ of the current weighted least squares problem.  To manage this difficulty, we shall consider a very specific recipe for generating the weights.  Other recipes are certainly possible.

Given a real 
number $\epsilon>0$ and a weight vector $w \in \R^N$, with $w_j>0$, $j=1,\dots,N$,  we define
\beqn
\label{defJ}
\cJ(z,w,\epsilon):= \frac{1}{2} \left [ \sum_{j=1}^N z_j^2w_j
+\sum_{j=1}^N(\epsilon^2w_j+w_j^{-1})\right ],\quad z\in \R^N.
\eeqn
Given $w$ and $\epsilon$, the element $z\in\R^N$   that minimizes 
$\cJ$ is unique because $\cJ$ is strictly convex.

Our algorithm will use an alternating method for choosing minimizers and weights based on the functional $\cJ$.  To describe this, we define for   $z\in\R^N$ the   non-increasing 
rearrangement $r(z)$  of the absolute values of the entries of $z$.
Thus $r(z)_i$ is the $i$-th largest element of the set 
$\{|z_j|, ~j=1,\dots,N\}$, and a vector $v$ is $k$-sparse iff $r(v)_{k+1}=0$.  

\begin{alg}
\label{alg1} 
{\rm We initialize by taking $w^0:=(1,\dots,1)$.  We also set $\epsilon_0:=1$.
We then recursively define for $n = 0,1,\dots,$
\beqn
\label{xn}
x^{n+1}:=\argmin_{z\in\cF(y)} ~\cJ(z,w^n,\epsilon_n)
= \argmin_{z\in\cF(y)} \|z\|_{\ell_2(w^n)}
\eeqn
and
\beqn
\label{en}
\epsilon_{n+1}:=\min (\epsilon_n, \frac{r(x^{n+1})_{K+1}}{N}),
\eeqn
where $K$ is a fixed integer that will be described more fully later.
We also define
\beqn
\label{wn}
w^{n+1}:=\argmin_{w>0} ~\cJ(x^{n+1},w,\epsilon_{n+1}).
\eeqn 
We stop the  algorithm if $\epsilon_n=0$; in this case we define $x^j:=x^n$ for 
$j>n$.  However, in general, the algorithm will generate an infinite sequence 
$(x^n)_{n\in \N}$ of distinct vectors. }\hfill $\square$
\end{alg}

Each step of the algorithm requires the solution of a  least squares problem.  In matrix form
\beqn
\label{mfsol}
x^{n+1}= D_n\Phi^t(\Phi D_n\Phi^t)^{-1}y,
\eeqn
where $D_n$ is the $N\times N$ diagonal matrix whose $j$-th diagonal entry is $w_j^n$ and $A^t$ denotes the transpose of the matrix $A$.  
Once $x^{n+1}$ is found, the weight $w^{n+1}$ is given by
\beqn
\label{wn1}
w_j^{n+1}=[(x_j^{n+1})^2+\epsilon_{n+1}^2]^{-1/2},\quad j=1,\dots,N.
\eeqn

 We shall prove several results about the convergence and rate of convergence of this algorithm.  This will be done under the following assumption on $\Phi$.
 
{\bf The Restricted Isometry Property (RIP):}  We say that the matrix $\Phi$ satisfies the Restricted Isometry Property of order $L$ with constant $\delta\in(0,1)$ if for each vector $z$ with sparsity $L$ we have
 \beqn
 \label{RIP}
 (1-\delta)\|z\|_{\ell_2^N}\leqslant \|\Phi z\|_{\ell_2^m}\leqslant (1+\delta)\|z\|_{\ell_2^N}.
 \eeqn
The RIP was introduced by Cand\`es and Tao \cite{cata05,cataXX} in their study of compressed sensing and $\ell_1$-minimization.  It has several analytical and geometrical interpretations that will be discussed in \S \ref{ripsect}.  To mention just one of these results (see \cite{CDDXX}), it is known that if $\Phi$ has the RIP of order $L:=J+J'$, with     
$\delta <\frac{\sqrt{J'}-\sqrt{J}}{\sqrt{J'}+\sqrt{J}}$ (here $J'>J$) and if 
\eref{system} has a $J$-sparse solution $z\in\cF(y)$, then 
this solution is the unique $\ell_1$ minimizer in $\cF(y)$. 
(This can still be sharpened: in \cite{caIMA}, Cand\`es showed that if
$\cF(y)$ contains a $J$-sparse vector, and 
if $\Phi$ has RIP of order $2J$ with $\delta<\sqrt{2}-1$,
then that $J$-sparse vector is unique and is the unique $\ell_1$ minimizer
in $\cF(y)$.)

The main result of this paper (Theorem \ref{conv}) is that whenever $\Phi$ satisfies the RIP of order $K+K'$ (for some $K'>K$) and $\delta$ sufficiently close 
to zero, then Algorithm \ref{alg1}
converges to a solution $\bar x$ of \eref{system} for each $y\in \R^m$. 
Moreover, if there is a  solution $z$ to \eref{system} that has sparsity $k\leqslant K-\kappa$, 
then  $\bar x=z$.  Here $\kappa>1$ depends on the RIP constant $\delta$ and can be made arbitrarily close to $1$ when $\delta$ is made small. 
The result cited in our previous paragraph implies that
in this case
$\bar x=x^*$, where $x^*$ is the $\ell_1$-minimal solution to \eref{system}.

A second part of our analysis concerns rates of convergence.  
We shall show that if \eref{system} has a $k$-sparse solution with, e.g., 
$k\leqslant K-4$ and if $\Phi$ satisfies the RIP of order $3K$ with 
$\delta$ sufficiently close 
to zero, then Algorithm 1
converges exponentially fast to $\bar x=x^*$.  Namely, once $x^{n_0}$ is sufficiently close to its limit $\bar x$, we have
\beqn
\label{exponential}
\|\bar x-x^{n+1}\|_{\ell_1^N}\leqslant \mu \|\bar x-x^n\|_{\ell_1^N},\quad n\geqslant n_0,
\eeqn
where $\mu<1$ is a fixed constant (depending on $\delta$). From this result it follows that we have exponential convergence to $\bar x$ whenever $\bar x$ is $k$-sparse; however we have no real information on how long
it will take before the iterates enter the region where we can control 
$\mu$. (Note that this is similar to convergence results for  the interior
point algorithms that can be used for direct $\ell_1$-minimization.)

The potential of IRLS algorithms, tailored to mimic $\ell_1$-minimization and so recover sparse solutions, 
has recently been investigated numerically by Chartrand and several
co-authors \cite{ch07,ch08,chyi08}. Our work provides proofs 
of several findings listed in these works.

One of the virtues of our approach is that, with minor technical modifications, 
it allows a similar detailed analysis of IRLS algorithms with weights that promote the {\em non-convex} optimization of $\ell_\tau$-norms for $0<\tau<1$. 
We can show not only that these algorithms can again recover sparse solutions, 
but also that their local rate of convergence is superlinear and tends 
to quadratic when $\tau$ tends to zero. Thus we also justify theoretically the recent numerical results by Chartrand et al.
concerning such non-convex  
$\ell_\tau$-norm optimization \cite{ch07,ch08,chst08,sachyi08}.

An outline of our paper is the following.  
In the next section we make some remarks about $\ell_1$- and weighted 
$\ell_2$-minimization, upon which we shall call in our proof.  
In the following section, we recall 
the Restricted Isometry Property and the Null Space Property including some of 
its consequences that are important to our analysis.  
In section \ref{prelim}, we gather 
some preliminary results we shall need to prove our main convergence 
result, Theorem \ref{conv}, 
which is formulated and proved
in section \ref{sectmain}.
We then turn to the issue on rate of convergence in 
section \ref{sectrate}.
In section \ref{sectnonconv} we 
generalize the convergence results obtained 
for $\ell_1$-minimization to the case of $\ell_\tau$-spaces for $0<\tau<1$; in 
particular, we show, with Theorem \ref{superth}, the local superlinear convergence 
of the IRLS algorithm in this setting. We conclude the paper with a short section 
dedicated to a few numerical examples that dovetail nicely with the theoretical results.

\section{Characterization of $\ell_1$- and weighted $\ell_2$-minimizers}
\label{unsection}
We fix $y\in\R^m$ and consider the underdetermined system $\Phi x=y$.   Given  a norm $\|\cdot \|$,  the problem of minimizing $\|z\|$ over $z\in \cF(y)$ can be viewed as a problem of approximation.  Namely, for any $x_0\in\cF(y)$, we can
characterize the minimizers in $\cF(y)$ as exactly those elements 
$z\in\cF(y)$ that can be written as  $z=x_0+\eta$, with $\eta$ a best approximation to $ - x_0$  from $\cN$. 
In this way one can characterize minimizers $z$ from classical results on best approximation in normed spaces.  We consider two examples of this in the present section, corresponding to  the $\ell_1$-norm and the weighted $\ell_2(w)$-norm.

Throughout this paper, we shall denote by $x$ any element from $\cF(y)$ that has smallest $\ell_1$-norm, as in \eref{l1}.  When $x$ is unique, we shall emphasize this by denoting it by $x^*$. In general, $x$ and $x^*$ need not be sparse, although
we will often consider cases where they are.
We begin with the following well-known lemma (see for example Pinkus \cite{pi89}) which characterizes the minimal $\ell_1$-norm elements from $\cF(y)$.
\begin{lemma}
\label{l1lemma}
An element $x\in \cF(y)$ has minimal $\ell_1$-norm among all elements $z\in\cF(y)$ if and only if 
\beqn
\label{l1char}
|\sum_{x_i\neq 0}\sign(x_i)\eta_i|\leqslant \sum_{x_i=0}|\eta_i|,\quad \eta\in\cN.
\eeqn
Moreover, $x$ is unique if and only if we have strict inequality in {\rm \eref{l1char}} for all $\eta\in\cN$ which are not identically zero.
\end{lemma}
\begin{Proof}  We give the simple proof for completeness of this paper.
If  $x\in\cF(y)$ has minimum $\ell_1$-norm, then we have, for any $\eta\in \cN$ and any $t\in\R$, 
\beqn
\label{l1char1}
\sum_{i=1}^N|x_i+t\eta_i|\geqslant \sum_{i=1}^N|x_i|.
\eeqn
Fix $\eta\in\cN$. If $t$ is sufficiently small then $x_i+t\eta_i$ and $x_i$ will have the same sign $s_i:=\sign (x_i)$ whenever $x_i\neq 0$. 
 Hence, \eref{l1char1} can be written as
$$
t\sum_{x_i\neq 0}s_i\eta_i+ \sum_{x_i=0}|t\eta_i|\geqslant 0.
$$
Choosing $t$ of an appropriate sign,  we see that \eref{l1char} is a necessary condition.

For the opposite direction, we note that if \eref{l1char} holds then for each $\eta\in\cN$, we have
\begin{eqnarray}
\label{l1char2}
\sum_{i=1}^N|x_i|&=& \sum_{x_i\neq 0}s_ix_i= \sum_{x_i\neq 0} 
s_i(x_i+\eta_i)-\sum_{x_i\neq 0}s_i\eta_i\cr
&\leqslant &\sum_{x_i\neq 0} s_i(x_i+\eta_i)+\sum_{x_i=0} |\eta_i|\leqslant  
\sum_{i=1}^N|x_i+\eta_i|, 
\end{eqnarray}
where the   first inequality uses \eref{l1char}.

If $x$ is unique then we have strict inequality in \eref{l1char1} and hence subsequently in \eref{l1char}.  If we have strict inequality in
\eref{l1char} then the subsequent strict inequality in \eref{l1char2}
implies uniqueness. 
\end{Proof}

\begin{remark}
\label{vanishing}
{\rm Applying 
Lemma \ref{l1lemma} to the special case of $\ell_1$-minimizers with no vanishing entries, we see that a vector
$x \in \mathcal F(y)$, with  $x_i \neq 0$ for all $i=1,\dots,N$, is a minimal 
$\ell_1$-norm solution if and only if 
\beqn
\label{zerochar}
\sum_{i=1}^N s_i \eta_i = 0, \quad \mbox{for all } \eta \in \mathcal N.
\eeqn
This implies that a minimal $\ell_1$-norm solution to $\Phi x =y$ for which all entries are non-vanishing
is necessarily non-unique, by the following argument. Suppose that $x_i \neq 0$ for all $i=1,\dots,N$ and that $x \in \mathcal F(y)$ is a minimal $\ell_1$-norm solution. Pick now any $\eta \in \mathcal N$, $\eta \neq 0$, and pick $t>0$
so that $t < \min_{\eta_i \neq 0} |x_i|/|\eta_i| $; it then follows
that $s_i = \sign(x_i + t \eta_i)$ for all $i=1,\dots,N$. But then we have  $\sum_{i=1}^N |x_i + t \eta_i|= \sum_{i=1}^N s_i (x_i + t \eta_i) = \sum_{i=1}^N |x_i|$ by \eref{zerochar}, so that $x + t \eta$ is also a minimal solution, different from $x$. Hence, unique $\ell_1$-minimizers are necessarily $k$-sparse for some $k < N$.} \hfill $\square$
\end{remark}

We next consider minimization in a weighted $\ell_2(w)$-norm.
We suppose that the weight $w$ is {\em strictly positive} which we define
to mean that
$w_j>0$ for all $j\in \{1,\dots, N\}$. 
In this case, $\ell_2(w)$ is a Hilbert space with the inner product
\beqn
\label{ip}
\<u,v\>_w:=\sum_{j=1}^Nw_ju_jv_j.
\eeqn
We define
\beqn
\label{xw}
x^w:=\argmin_{z\in\cF(y)}\|z\|_{\ell_2^N(w)}.
\eeqn
Because the $\|\cdot\|_{\ell_2^N(w)}$-norm is strictly convex, the 
minimizer $x^w$ is necessarily unique; it is   completely characterized by the orthogonality conditions
\beqn
\label{wortho}
\langle x^w,\eta\rangle_w=0,\quad \forall \eta\in\cN.
\eeqn
Namely, $x^w$ necessarily satisfies \eref{wortho}; on the other hand, any element $z \in \cF(y)$ that satisfies $\langle z,\eta\rangle_w=0$ for all $\eta\in\cN$ is automatically equal to $x^w$.

At this point, we would like to tabulate some of the notation we have used
in this paper to denote various kinds of 
minimizers and other solutions alike (such as limits of algorithms).

\begin{table}[h]
\begin{center}
\begin{tabular}{|l|l|}
\hline
$z$ & an (arbitrary) element of $\cF(y)$ \\
\hline
$x$ & any solution of $\displaystyle \min_{z \in \cF(y)} \|z\|_{\ell_1}$ \\
\hline
$x^*$ & unique solution of $\displaystyle \min_{z \in \cF(y)} \|z\|_{\ell_1}$ 
(notation used only when the minimizer is unique) \\
\hline
$x^w$ & unique solution of $\displaystyle \min_{z \in \cF(y)} \|z\|_{\ell_2(w)}$, $w_j > 0$ for all $j$ \\
\hline
$\bar{x}$ & limit of Algorithm \ref{alg1}\\
\hline
$x^\e$ & unique solution of $\displaystyle \min_{z \in \cF(y)} f_\e(z)$; 
see \eref{newfunct}\\
\hline
\end{tabular}
\caption{Notation for solutions and minimizers.}
\end{center}
\end{table}

\section{The Restricted Isometry and the Null Space Properties}
\label{ripsect}
To analyze the convergence of our algorithm, we shall impose the Restricted Isometry Property (RIP) already mentioned in the introduction, or a slightly
weaker version, the {\em Null Space Property}, which will be defined below.     
Recall that  $\Phi$ satisfies RIP of order $L$ for $\delta\in(0,1)$ (see \eref{RIP})
iff
\beqn
\label{RIP1}
(1-\delta) \|z\|_{\ell_2^N} \leqslant  \| \Phi z \|_{\ell_2^m} \leqslant (1+\delta) \|z\|_{\ell_2^N},\quad \mbox{ for all }L \mbox{-sparse } z.
\eeqn
It is known that many families of matrices satisfy the RIP.  While there are deterministic families that are known to satisfy RIP, the largest range of $L$, 
(asymptotically, as $N\to\infty$, with e.g. $m/N$ kept constant) is obtained (to date) by using random families.
For example, random families in which the entries of the matrix $\Phi$ are independent realizations of a (fixed) Gaussian or Bernoulli random variable are known to have the RIP with high probability for each $L\leqslant c_0(\delta)\,n/\log n$ (see \cite{cata05,cataXX,BDDW,ruveXX} for a discussion of these results).

We shall say
that $\Phi$ has the 
{\em Null Space Property} (NSP) of order $L$ for $\gamma > 0$ if  \footnote{This definition of the Null Space Property is a slight variant of that given in \cite{CDDXX} but is more convenient for the results in the present paper.}
\beqn
\label{NSP}
\|\eta_T\|_{\ell_1}\leqslant \gamma\|\eta_{T^c}\|_{\ell_1},
\eeqn
for all sets $T$ of cardinality not exceeding $L$ and all $\eta \in \cN$.  
Here and later, we denote by $\eta_S$ the vector obtained from $\eta$ by 
setting to zero all coordinates $\eta_i$ for 
$i \notin S \subset \{1,2,\ldots,N\}$;  
$T^c$ denotes the complement of the set $T$.
It is shown in Lemma 4.1 of \cite{CDDXX} that if $\Phi$ 
has the RIP of order $L:=J+J'$ 
for a given $\delta\in(0,1)$, where $ J,J' \geqslant 1$ are integers, 
then $\Phi$ has the NSP of order $K$ for
$\gamma:=\frac{1+\delta}{1-\delta}\,\sqrt{\frac{J}{J'}}$.  
Note that if $J'$ is sufficiently large then $\gamma<1$. 

Another result in \cite{CDDXX} (see also Lemma \ref{NSPl1min} below)
states 
that in order to guarantee that a $k$-sparse
vector $x^*$ is the unique $\ell_1$-minimizer in $\cF(y)$, it is
sufficient that $\Phi$ has the NSP of order $L \geqslant  k$ and $\gamma<1$.
(In fact, the argument in \cite{cataXX}, proving that for $\Phi$
with the RIP, $\ell_1$-minimization
identifies sparse vectors in $\cF(y)$, can be split into two steps: one that
implicitly derives the NSP 
from the RIP, and the remainder of the proof, which uses only the NSP.) \\ 

 Note that if the NSP holds for some order $L_0$ and constant $\gamma_0$ (not necessarily $<1$),
then, by choosing
$a>0$ sufficiently small, one can ensure that $\Phi$ has the NSP of order 
$L = a L_0$ with constant 
$\gamma<1$ (see \cite{CDDXX} for details).
So the effect of requiring that $\gamma<1$ is tantamount to reducing the range of 
$L$ slightly.
\\

When proving results on the convergence of our algorithm later in this paper, we shall state them under the assumptions that $\Phi$ has the NSP for some $\gamma<1$ and an appropriate value of $L$. Using the observations above, they can easily be rephrased in terms of RIP bounds for $\Phi$.

\section{Preliminary results}
\label{prelim}
We first make some comments about the decreasing rearrangement $r(z)$ and
the $j$-term approximation errors for vectors in $\R^N$. 
Let us denote by $\Sigma_k$ the set of all $x\in \R^N$ such that
$\#(\mbox{supp}(x))\leqslant k$.  
For any 
$z \in \R^N$ and any $j=1,2,\dots, N$, we denote by
\beqn
\label{sigma}
\sigma_j(z)_{\ell_1}:= \inf_{w\in\Sigma_j}\|z-w\|_{\ell_1^N}
\eeqn
the $\ell_1$-error in approximating a general vector $z\in\R^N$ by a $j$-sparse vector. Note that these approximation errors can be written as a sum of
entries of $r(u)$:  $\sigma_j(z)_{\ell_1}=\sum_{ \nu>j}r(z)_\nu$.
We have the following lemma:

\begin{lemma}
\label{dlemma}
The map $z \mapsto r(z)$ is Lipschitz continuous on $(\R^N, \|\cdot
\|_{\ell_\infty})$: for any $z,z'\in\R^N$, we have
\beqn
\label{dl1}
\|r(z)-r(z')\|_{\ell_\infty} \leqslant 
\|z-z'\|_{\ell_\infty}.
\eeqn
Moreover, for any $j$, we have
\beqn
\label{dl10}
|\sigma_j(z)_{\ell_1}-\sigma_j(z')_{\ell_1}| \leqslant \|z-z'\|_{\ell_1},
\eeqn
and for any $J > j$, we have
\beqn
\label{dl11}
(J-j)r(z)_{J}\leqslant \|z-z'\|_{\ell_1}+\sigma_{j}(z')_{\ell_1}.
\eeqn
\end{lemma}

\begin{Proof}
For any pair of points $z$ and $z'$, and 
any $j \in \{1,\dots,N\}$, let $\Lambda$ be a set of $j-1$ indices
corresponding to the $j-1$ largest entries in $z'$. 
Then
\beqn
\label{dl2}
r(z)_j \leqslant  \max_{i\in \Lambda^c}|z_i|\leqslant \max_{i\in \Lambda^c}|z_i'| 
+\|z-z'\|_{\ell_\infty} = r(z')_j +\|z-z'\|_{\ell_\infty}.
\eeqn
We can also reverse the roles of $z$ and $z'$.  Therefore,  we obtain \eref{dl1}.
To prove \eref{dl10}, we approximate $z$ by a $j$-term best approximation 
$u\in \Sigma_j$ of $z'$ in $\ell_1$. Then 
\beqns
\sigma_j(z)_{\ell_1} 
\leqslant \|z-u\|_{\ell_1} \leqslant \|z-z'\|_{\ell_1} + \sigma_j(z')_{\ell_1},
\eeqns
and the result follows from symmetry.

To prove \eref{dl11}, it   suffices to note that 
$(J-j)\,r(z)_{J} \leqslant \sigma_j(z)_{\ell_1}$.
\end{Proof}

Our next result  is an approximate reverse triangle inequality for 
points in $\cF(y)$.  Its importance to us lies in its implication that whenever
two points $z,z'\in\cF(y)$  have close $\ell_1$-norms
and one of them is close to a $k$-sparse vector,
then they necessarily are close to each other. (Note that it also implies  
that the other vector must then also be close to that $k$-sparse vector.)
This is a geometric property 
of the null space. 
\begin{lemma}
\label{balllemma}
Assume that {\rm \eref{NSP}} holds for some $L$ and $\gamma<1$.  Then, for 
any  $z,z'\in\cF(y)$, we have
\beqn
\label{ball1}
\|z'-z\|_{\ell_1}\leqslant \frac{1+\gamma}{1-\gamma}
\left (\|z'\|_{\ell_1}-\|z\|_{\ell_1}+2\sigma_L(z)_{\ell_1}\right ).
\eeqn
\end{lemma}
 
\begin{Proof}  Let $T$ be a set of indices of the $L$ largest entries in 
$z$. Then
\begin{eqnarray}
\label{ball2}
\|(z'-z)_{T^c}\|_{\ell_1}
&\leqslant& \|z'_{T^c}\|_{\ell_1}+\|z_{T^c}\|_{\ell_1} \cr
&=& \|z'\|_{\ell_1}-\|z'_T\|_{\ell_1}+\sigma_L(z)_{\ell_1}\cr
&=& \|z\|_{\ell_1}+\|z'\|_{\ell_1}-\|z\|_{\ell_1}-\|z'_T\|_{\ell_1}
+\sigma_L(z)_{\ell_1}\cr
&= & \|z_T\|_{\ell_1}-\|z'_T\|_{\ell_1}+\|z'\|_{\ell_1}-\|z\|_{\ell_1}+2
\sigma_L(z)_{\ell_1}\cr
&\leqslant& \|(z'-z)_T\|_{\ell_1}+\|z'\|_{\ell_1}-\|z\|_{\ell_1}+2\sigma_L(z)_{\ell_1}.
\end{eqnarray}
Using \eref{NSP}, this gives 
\beqn
\label{ball3}
\|(z'-z)_T\|_{\ell_1} 
\leqslant \gamma  \|(z'-z)_{T^c}\|_{\ell_1}
\leqslant \gamma ( \|(z'-z)_T\|_{\ell_1}+\|z'\|_{\ell_1}-\|z\|_{\ell_1}
+2\sigma_L(z)_{\ell_1}).
\eeqn
In other words,
\beqn
\label{ball4}
\|(z'-z)_T\|_{\ell_1}\leqslant  \frac{\gamma}{1-\gamma} 
(\|z'\|_{\ell_1}-\|z\|_{\ell_1}+2\sigma_L(z)_{\ell_1}).
\eeqn
Using this, together with \eref{ball2}, we obtain
\beqn
\label{ball5}
\|z'-z\|_{\ell_1} = \|(z'-z)_{T^c}\|_{\ell_1} +\|(z'-z)_T\|_{\ell_1}
\leqslant \frac{1+\gamma}{1-\gamma} 
( \|z'\|_{\ell_1}-\|z\|_{\ell_1}+2\sigma_L(z)_{\ell_1}),
\eeqn
as desired.  
\end{Proof}

This result then allows the following simple proof
of some of the results of \cite{CDDXX}:
\begin{lemma}
\label{NSPl1min}
Assume that {\rm \eref{NSP}} holds for some $L$ and $\gamma<1$. Suppose
that $\cF(y)$ contains an $L$-sparse vector. Then this vector is the
unique $\ell_1$-minimizer in $\cF(y)$; denoting it by $x^*$, we have
moreover, for all $v\in \cF(y)$, 
\beqn
\label{ball17}
\|v-x^*\|_{\ell_1} 
\leqslant 2 \,\frac{1+\gamma}{1-\gamma} \,
\sigma_L(v)_{\ell_1}\,.
\eeqn
\end{lemma}
\begin{Proof}
For the time being, we denote the $L$-sparse vector in $\cF(y)$ by $x_s$.\\
Applying \eref{ball1} with $z'=v$ and $z=x_s$, we find
$$
\|v-x_s\|_{\ell_1}\leqslant  \frac{1+\gamma}{1-\gamma} [\|v\|_{\ell_1}-\|x_s\|_{\ell_1}] \, ;
$$
since $v\in \cF(y)$ is arbitrary, this implies that 
$\|v\|_{\ell_1}-\|x_s\|_{\ell_1}\geqslant 0$ for all $v\in\cF(y)$, so that
$x_s$ is an $\ell_1$-norm
minimizer in $\cF(y)$.

If $x'$ were another $\ell_1$-minimizer in $\cF(y)$, then it would follow that
$\|x'\|_{\ell_1}=\|x_s\|_{\ell_1}$, and the inequality
we just derived would imply $\|x'-x_s\|_{\ell_1}=0$, or $x'=x_s$. It 
follows that $x_s$ is the unique $\ell_1$-minimizer in $\cF(y)$, which
we denote by $x^*$, as proposed earlier.

Finally, we apply \eref{ball1} with $z'=x^*$ and $z=v$, and we obtain
$$
\|v-x^*\|\leqslant  \frac{1+\gamma}{1-\gamma} (\|x^*\|_{\ell_1}-\|v\|_{\ell_1}+2 \sigma_L(v)_{\ell_1}) 
\leqslant  2 \frac{1+\gamma}{1-\gamma} \sigma_L(v)_{\ell_1}\, ,
$$
where we have used the $\ell_1$-minimization property of $x^*$. 
\end{Proof}

Our next set of remarks centers around the functional $\cJ$ defined by 
\eref{defJ}.  Note that for each $n=1,2,\dots$, we have
\beqn
\label{eqvJ}
\cJ(x^{n+1},w^{n+1},\epsilon_{n+1}) = \sum_{j=1}^N [(x_j^{n+1})^2
+ \epsilon_{n+1}^2]^{1/2}.
\eeqn
We also have the following monotonicity property which holds for 
all $n\geqslant 0$:
\beqn
\label{mono}
\cJ(x^{n+1},w^{n+1},\epsilon_{n+1})\leqslant \cJ(x^{n+1},w^{n},\epsilon_{n+1})
\leqslant\cJ(x^{n+1},w^{n},\epsilon_{n})\leqslant\cJ(x^{n},w^{n},\epsilon_{n}).
\eeqn
Here the first inequality follows from 
the minimization property that defines $w^{n+1}$, 
the second inequality from $\epsilon_{n+1}\leqslant \epsilon_n$, 
and the last inequality from 
the minimization property that defines $x^{n+1}$. 
For each $n$, $x^{n+1}$ is completely determined by
$w^n$; for $n=0$, in particular, $x^1$ is determined solely by 
$w^0$, and independent of the choice of $x^0 \in \cF(y)$.
(With the initial weight vector defined by
$w^0=(1,\dots,1)$, 
$x^1$ is the classical minimum $\ell_2$-norm element of $\cF(y)$.)
The inequality \eref{mono} for $n=0$ thus holds for 
arbitrary $x^0 \in \cF(y)$.
\begin{lemma} 
\label{lemmawb} For each $n \geqslant 1$ we have
\beqn
\label{xnbound}
\|x^n\|_{\ell_1}\leqslant \cJ(x^1,w^0,\epsilon_0)=:A
\eeqn
and
\beqn
\label{weightbound}
w_j^n\geqslant A^{-1}, \quad j=1,\dots,N.
\eeqn
\end{lemma}

\begin{Proof}  
The bound \eref{xnbound} follows from \eref{mono} and   
$$
\|x^n\|_{\ell_1}\leqslant \sum_{j=1}^N[(x_j^n)^2+\epsilon_n^2]^{1/2}= \cJ(x^n,w^n,
\epsilon_n).
$$
The bound
\eref{weightbound} follows from 
 $(w_j^n)^{-1}=[(x_j^n)^2+\epsilon_n^2]^{1/2}
\leqslant \cJ(x^n,w^n,\epsilon_n)\leqslant   A$, 
where the last inequality uses \eref{mono}.
\end{Proof}
 
\section{Convergence of the algorithm}
\label{sectmain}

In this section, we prove that the algorithm converges.  Our starting point is the following lemma 
that establishes  $(x^n-x^{n+1}) \to 0$ for $n \to \infty$.

\begin{lemma}
\label{boundlemma}
Given any $y\in \R^m$, the $x^n$ satisfy
\beqn
\label{bl1}
\sum_{n=1}^\infty \|x^{n+1}-x^n\|_{\ell_2}^2\leqslant 2 A^2.
\eeqn
where $A$ is the constant of Lemma {\rm \ref{lemmawb}}. In particular,
we have 
\beqn
\label{asympreg}
\lim_{n \to \infty} (x^n - x^{n+1}) = 0.
\eeqn
\end{lemma}

\begin{Proof} For each $n=1,2,\dots$, we have
\begin{eqnarray}
\label{bl2}
2[\cJ(x^n,w^n,\epsilon_n)-\cJ(x^{n+1},w^{n+1},\epsilon_{n+1})] 
&\geqslant&  
2[\cJ(x^n,w^n,\epsilon_n)-\cJ(x^{n+1},w^n,\epsilon_n)] \cr
&=& \langle x^n,x^n\rangle_{w^n}-\langle x^{n+1}, x^{n+1}\rangle_{w^n}\cr
&=&\langle x^n+x^{n+1},x^n-x^{n+1}\rangle_{w^n} \cr
&=&\langle x^n-x^{n+1}, x^n-x^{n+1}\rangle_{w^n}\cr
&=&\sum_{j=1}^Nw_j^n(x_j^n-x_j^{n+1})^2\cr
&\geqslant& A^{-1} \| x^n-x^{n+1}\|_{\ell_2}^2,
\end{eqnarray}
where the third equality uses the fact that
$\langle x^{n+1},x^{n}-x^{n+1}\rangle_{w^n}=0$ (observe that 
$x^{n+1}-x^{n} \in \cN$ and invoke \eref{wortho}), 
and the inequality uses the bound \eref{weightbound} on the weights.  If we 
now sum these inequalities over $n\geqslant 1$, we arrive
at \eref{bl1}.
\end{Proof}

From the monotonicity of $\epsilon_n$, we know that 
$\epsilon:=\lim_{n\to \infty}\epsilon_n$ exists and is non-negative. 
The following functional will play an important role in our proof of 
convergence:
\beqn
\label{newfunct}
f_\epsilon(z):=\sum_{j=1}^N(z_j^2+\epsilon^2)^{1/2}.
\eeqn
Notice that if we knew that $x^n$ converged to $x$ then, in view of 
\eref{eqvJ}, $f_\epsilon(x)$ would be the limit of 
$\cJ(x^{n},w^{n},\epsilon_n)$.   
When $\epsilon > 0$ the functional $f_\epsilon$  is strictly 
convex and therefore has a unique minimizer 
\beqn
\label{xe}
x^\e:=\argmin_{z\in\cF(y)}f_\epsilon(z).
\eeqn 
This minimizer is characterized by the following lemma:
\begin{lemma}
\label{min_f_eps}
Let $\e > 0$ and $z \in \cF(y)$. Then
$z=x^\e$ if and only if $\langle z, \eta
\rangle_{\widetilde w(z,\e)} = 0$ for all $\eta \in \cN$, 
where $\widetilde w(z,\e)_i = 
[z_i^2 + \e^2]^{-1/2}$.
\end{lemma}
\begin{Proof}
For the ``only if'' part, let $z = x^\e$ and
$\eta \in \cN$ be arbitrary. Consider the analytic function
\beqns 
G_\e(t) := f_\e(z + t \eta) - f_\e(z).
\eeqns
We have $G_\e(0) = 0$, and 
by the minimization property $G_\e(t) \geqslant 0$ for all $t \in \R$.
Hence,  $G'_\e(0) = 0$. A simple calculation reveals that
\beqns
G'_\e(0) = 
\sum_{j=1}^N \frac{\eta_i z_i}{[z_i^2 + \e^2]^{1/2}} 
= \langle z, \eta \rangle_{\widetilde w(z,\e)},
\eeqns 
which gives the desired result.

For the ``if'' part, assume that $z \in \cF(y)$ and
$\langle z, \eta
\rangle_{\widetilde w(z,\e)} = 0$ for all $\eta \in \cN$, where 
$\widetilde w(z,\e)$ is
defined as above.
We shall show that $z$ is a minimizer of $f_\epsilon$ on $\cF(y)$.  
Indeed, consider the convex univariate function  $[u^2+\epsilon^2]^{1/2}$. 
For any point $u_0$ we have from convexity that
\beqn
\label{convex}
[u^2+\epsilon^2]^{1/2}\geqslant  [u_0^2+\epsilon^2]^{1/2}+[u_0^2+\epsilon^2]^{-1/2}
u_0(u-u_0),
\eeqn
because the right side is the linear function which is tangent to this 
function at $u_0$. 
It follows that for any point $v\in\cF(y)$ we have
\beqn
\label{first}
f_\epsilon(v)\geqslant f_\epsilon(z) + 
\sum_{j=1}^N [z_j^2+\epsilon^2]^{-1/2} z_j(v_j-z_j)
= f_\epsilon(z) + 
\langle z, v-z\rangle_{\tilde w(z,\e)} 
= f_\epsilon(z),
\eeqn
where we have used the orthogonality condition \eref{ortho2} and the fact 
that  $v-z$ is in $\cN$. 
Since $v$ is arbitrary, it follows that $z=x^\e$, as claimed.
\end{Proof}

We now give the  convergence of the algorithm.
\\

\begin{theorem} \label{conv} 
Let $K$ (the same index as used in the update rule {\rm \eref{en}}) be chosen so that  $\Phi$ satisfies the Null
Space Property {\rm \eref{NSP}} of order $K$, with $\gamma<1$.
Then, for each $y\in\R^m$, the output of Algorithm \ref{alg1} converges 
to a vector $\bar x$, with $r(\bar x)_{K+1}= N \lim_{n \to \infty} \e_n$ 
and  the following hold:\\
{\rm (i)} If $\e=\lim_{n \to \infty}\e_n =0$, then $\bar x$ is $K$-sparse; in this
case there is therefore a unique $\ell_1$-minimizer $x^*$, and $\bar x = x^*$; 
moreover, we have, for $k \leqslant K$, and any $z \in \cF(y)$,
\beqn
\label{estt}
\|z-\bar x\|_{\ell_1}\leqslant c\sigma_k(z)_{\ell_1},\quad \mbox{ with }
c:=\frac{2(1+\gamma)}{1-\gamma}
\eeqn
{\rm (ii)} If $\e=\lim_{n \to \infty}\e_n >0$, then $\bar x = x^{\e}$;\\
{\rm (iii)} In this last case, if $\gamma$ satisfies the stricter bound 
$\gamma < 1-\frac{2}{K+2}$ (or,
equivalently, if $\frac{2\gamma}{1-\gamma}<K$), then we have, for all 
$z \in \cF(y)$ and any  $k< K-\frac{2\gamma}{1-\gamma}$, that
\beqn
\label{est}
\|z-\bar x\|_{\ell_1}\leqslant \tilde c \sigma_k(z)_{\ell_1},\quad \mbox{ with }
\tilde c :=\frac{2(1+\gamma)}{1-\gamma}\left[
\frac{K-k+\frac{3}{2}}{K-k-\frac{2\gamma}{1-\gamma}}\right]
\eeqn
As a consequence, this case is excluded  if $\cF(y)$ contains a vector of sparsity 
$k< K-\frac{2\gamma}{1-\gamma}$.
\end{theorem}
\noindent
The constant $\tilde c$ can be quite reasonable; 
for instance, if $\gamma \leqslant 1/2$ and $k\leqslant K-3$, then
we have $\tilde c \leqslant 9\, \frac{1+\gamma}{1-\gamma} \leqslant 27$.\\

\begin{Proof}
Note that since $\e_{n+1}\leq\e_n$, the $\e_n$ always converge. 
We start by considering the case $\epsilon:=\lim_{n\rightarrow \infty}
\epsilon_n=0$. 

{\bf Case $\epsilon=0$:} In this case, we want to prove that
$x^n$ converges , and that its limit is an $\ell_1$-minimizer. 
Suppose that $\epsilon_{n_{{\!\,}_0}} = 0$ for some $n_0$. 
Then by the definition of the algorithm, we know that 
the iteration is stopped at $n=n_0$, and 
$x^n=x^{n_{{\!\,}_0}}$, $n\geqslant n_0$. Therefore $\bar x = x^{n_0}$.
From the definition of $\e_n$, it then also follows that 
$r(x^{n_{{\!\,}_0}})_{K+1}=0$ and so $\bar x=x^{n_{{\!\,}_0}}$ is $K$-sparse.  As noted in \S\ref{ripsect} and Lemma \ref{NSPl1min}, if a $K$-sparse solution
exists when $\Phi$ satisfies the NSP of order $K$ with $\gamma < 1$,
then it is the unique $\ell_1$-minimizer.
Therefore, $\bar x$ equals $x^*$, this unique minimizer.

Suppose now  that $\epsilon_n > 0$ for all $n$.
Since $\e_n \to 0$, 
there is an increasing sequence of indices $(n_i)$ such that
$\e_{n_i} < \epsilon_{n_i-1}$ for all $i$. By the definition \eref{en}
of $(\e_n)_{n\in \N}$, we must
have $r(x^{n_i})_{K+1} < N \epsilon_{n_i-1}$ for all $i$.
Noting that $(x^{n})_{n\in\N}$ is a bounded sequence, there exists a subsequence
$(p_j)_{j \in \N}$ of $(n_i)_{i \in \N}$ such that $(x^{p_j})_{j\in \N}$ converges 
to a point
$\widetilde x \in \cF(y)$. By Lemma \ref{dlemma}, we know that 
$r(x^{p_j})_{K+1}$ converges to $r(\widetilde x)_{K+1}$. Hence we get
\beqn
\label{est4}
r(\widetilde x)_{K+1} = \lim_{j\to\infty}r(x^{p_j})_{K+1} 
\leqslant \lim_{j \to\infty} N\epsilon_{p_j-1} = 0, 
\eeqn
which means that the support-width of $\widetilde x$ is at most $K$, i.e.
$\widetilde x$ is $K$-sparse. 
By the same token used above, we again have that
$\widetilde x=x^*$, the unique $\ell_1$-minimizer.
We must still show that
$x^n \to x^*$. Since $x^{p_j} \to x^*$ and $\epsilon_{p_j} \to 0$, 
\eref{eqvJ} implies $\cJ(x^{p_j},w^{p_j},\epsilon_{p_j}) \to 
\|x^*\|_{\ell_1}$. By the monotonicity property stated in \eref{mono},
we get $\cJ(x^n,w^n,\e_n) \to \|x^*\|_{\ell_1}$. Since \eref{eqvJ} implies
\beqn
\cJ(x^n,w^n,\e_n) - N\e_n \leqslant \|x^n\|_{\ell_1} \leqslant  \cJ(x^n,w^n,\e_n),
\eeqn
we obtain $\| x^n \|_{\ell_1} \to \|x^*\|_{\ell_1}$.
Finally, we invoke Lemma \ref{balllemma} with $z'=x^n$, $z=x^*$, and $k=K$
to get
\beqn
\limsup_{n\to \infty} \|x^n-x^*\|_{\ell_1} \leqslant  \frac{1+\gamma}{1-\gamma} 
\left(\lim_{n \to \infty} \|x^n\|_{\ell_1} - \|x^*\|_{\ell_1} \right)
= 0,
\eeqn
which completes the proof that $x^n \to x^*$ in this case. 

Finally, \eref{estt} follows from \eref{ball17} of
Lemma \ref{NSPl1min} (with $L=K$), and the observation that 
$\sigma_n(z) \geqslant \sigma_{n'}(z)$ if $n \leqslant n'$.

{\bf Case $\epsilon>0$:} 
We shall first show that $x^n\to   x^\e$,  
$n\to \infty$, with $x^\e$ as defined by \eref{xe}. By Lemma \ref{lemmawb},
we know that  $(x^n)_{n=1}^\infty$   is  a bounded sequence in $\R^N$ and 
hence this sequence has accumulation points.  Let $(x^{n_i})$ be any 
convergent subsequence of $(x^n)$ and let $\widetilde x \in \cF(y)$ 
be its limit. We want to 
show that $\widetilde x = x^\e$.  

Since $w_j^n= [(x_j^n)^2+\epsilon_n^2]^{-1/2}\leqslant \epsilon^{-1}$, it follows 
that $\lim_{i\to\infty}w_j^{n_i} = [(\widetilde x_j)^2+\epsilon^2]^{-1/2}
=\widetilde w(\widetilde x, \e)_j$ $ =:\widetilde w_j$, $j=1,\dots,N$. On the other hand, by
invoking Lemma \ref{boundlemma}, we now 
find that $x^{n_i+1}\to \widetilde x$, $i \to\infty$.  
It then follows from the orthogonality relations \eref{wortho} that
for every $\eta\in \cN$, we have
\beqn
\label{ortho2}
\langle \widetilde x,\eta\rangle_{\widetilde w}  =
\lim_{i\to\infty}  \langle x^{n_i+1},\eta\rangle_{w^{n_i}} =0.
\eeqn
Now the ``if'' part of Lemma \ref{min_f_eps} implies that 
$\widetilde x = x^\e$. Hence $x^\e$ is the unique
accumulation point of $(x^n)_{n\in \N}$ and therefore its limit. This 
establishes (ii).

To prove the error estimate \eref{est} stated in (iii), we first note that for 
any $z \in \cF(y)$, we have
\beqn
\label{fnote}
\|x^\e\|_{\ell_1}\leqslant f_\epsilon(x^\e)
\leqslant f_\epsilon(z) \leqslant \|z\|_{\ell_1}+N\epsilon,
\eeqn
where the second inequality uses the minimizing property of $x^\e$.
Hence it follows that $\|x^\e\|_{\ell_1}-\|z\|_{\ell_1}\leqslant   N\epsilon$.
We now invoke Lemma \ref{balllemma} to obtain

\beqn
\label{c11}
\|x^\e-z\|_{\ell_1}\leqslant \frac{1+\gamma}{1-\gamma}
[N\epsilon+2\sigma_k(z)_{\ell_1}].
\eeqn
From Lemma \ref{dlemma} and \eref{en}, we obtain
\beqn
\label{eps}
N\epsilon
= \lim_{n\to \infty} N \epsilon_n
\leqslant \lim_{n\to\infty} r(x^n)_{K+1}=r(x^\e)_{K+1}.
\eeqn
It follows from \eref{dl11} that
\begin{eqnarray}
\label{est1}
(K+1-k)N\epsilon
&\leqslant &(K+1-k)r(x^\e)_{K+1} \cr
&\leqslant & \|x^\e-z\|_{\ell_1}+\sigma_k(z)_{\ell_1}\cr
&\leqslant &\frac{1+\gamma}{1-\gamma}
[N\epsilon+2\sigma_k(z)_{\ell_1}]+\sigma_k(z)_{\ell_1},
\end{eqnarray}
where the last inequality uses \eref{c11}.
Since by assumption on $K$, we have $K-k>\frac{2\gamma}{1-\gamma}$,
i.e. $K+1-k> 
\frac{1+\gamma}{1-\gamma}$, we obtain
\beqns
N \epsilon + 2\sigma_k(z)_{\ell_1}\leqslant 
\frac{2(K-k)+3}{(K-k)-\frac{2\gamma}{1-\gamma}} \, \sigma_k(z)_{\ell_1}.
\eeqns
Using this back in \eref{c11}, we arrive at \eref{est}.

Finally, notice that if $\cF(y)$ contains a $k$-sparse vector  (with 
$k < K -\frac{2\gamma}{1-\gamma}$), then we know already (see \S
\ref{ripsect})
that this must be the unique $\ell_1$-minimizer $x^*$; it
then follows from our arguments
above that we must have $\e=0$.
Indeed, if we had $\e>0$, then \eref{est1} would hold for $z=x^*$; since
$x^*$ is $k$-sparse, $\sigma_k(x^*)_{\ell_1}=0$, implying $\e=0$,
a contradiction with the assumption $\e>0$. This finishes the proof.

\end{Proof}

\begin{remark}
{\rm  Let us briefly compare our analysis of the IRLS algorithm with $\ell_1$ minimization.   The latter  recovers a $k$-sparse
solution (when one exists) if $\Phi$ has the NSP of order $K$ and $k\leqslant K$. 
  The analysis given in our proof of Theorem \ref{conv}
guarantees that our IRLS algorithm recovers $k$-sparse $x$ for a slightly smaller range of values $k$ than $\ell_1$-minimization, namely for 
$k <  K -\frac{2\gamma}{1-\gamma}$. Notice that this ``gap'' vanishes for vanishingly small $\gamma$. Although we have no examples to demonstrate, 
our arguments cannot exclude the case where
$\cF(y)$ contains 
a $k$-sparse vector $x^*$ with 
$K-\frac{2\gamma}{1-\gamma}\leqslant k \leqslant K$ (e.g., if
$\gamma \geqslant 1/3$ and $k=K-1$), and
our IRLS algorithm converges to $\bar x$, yet  $ \bar x \neq x^*$.  However, note that unless $\gamma$ is close to 1, the range of $k$-values in this ``gap" is fairly small; for instance, for $\gamma < \frac{1}{3}$, this
non-recovery of a $k$-sparse $x^*$ can happen only if $k =  K$. }\hfill 
$\square$
\end{remark}
\begin{remark}
\label{oneconstant}
{\rm The constant $c$ in \eref{estt} is clearly smaller than the constant $\tilde c$
in \eref{est}; it follows that when $k < K -\frac{2\gamma}{1-\gamma}$, the
estimate \eref{est} holds for all cases, regardless of whether $\e=0$ or not.}\hfill 
$\square$
\end{remark}
\section{Rate of Convergence}
\label{sectrate}
Under the conditions of Theorem \ref{conv}  the algorithm 
converges to a limit $\bar x$;  if there is a
$k$-sparse vector in $\cF(y)$ with $k< K-\frac{2\gamma}{1-\gamma}$, then this limit coincides with that
$k$-sparse vector, which is then also automatically the 
unique $\ell_1$-minimizer
$x^*$.
In this section our goal is to establish a bound for the 
rate of convergence in both the sparse and  non-sparse cases.  In the latter case, the goal is to establish the rate at which 
$x^n$ approaches to 
a ball of radius $C_1 \sigma_k(x^*)_{\ell^1}$ centered at
$x^*$.  We shall work under the same assumptions as in Theorem 
\ref{conv}.

\subsection{Case of $k$-sparse vectors}

Let us begin by assuming  that $\cF(y)$ contains the  $k$-sparse vector 
$x^*$.  The algorithm produces the sequence $x^n$, which converges to 
$x^*$, as established above.    Let us denote the (unknown) support 
of the 
$k$-sparse vector $x^* $ by $T$.

We introduce
an auxiliary sequence of error vectors
$\eta^n \in \cN$ via $\eta^n := x^n -  x^*$ and
\beqns
E_n := \| \eta^n \|_{\ell_1}=\|x^*-x^n\|_{\ell_1^N}. 
\eeqns
We know that $E_n \to 0$.   The following theorem gives a bound on the rate of convergence of $E_n$ to zero.

\begin{theorem} 
\label{k-sparse-rate}
 Assume $\Phi$ satisfies NSP of order $K$ with constant 
$\gamma$ such that $0<\gamma < 1 -\frac{2}{K+2}$.
Suppose that  
$k <  K -\frac{2\gamma}{1-\gamma}$ , $0<\rho<1$, and $0<\gamma< 1 -\frac{2}{K+2}$  are such that 
\beqns
\mu := \frac{\gamma(1+\gamma)}{1-\rho}\left(1  + \frac{1}{K+1-k}\right) < 1.
\eeqns
Assume that $\cF(y)$ contains a $k$-sparse vector $x^*$ and let $T=
\supp(x^*)$.
Let $n_0$ be such that
\beqn
\label{sball}
E_{n_{{\!\,}_0}} \leqslant   R^* := \rho \,\min_{i \in T} |x^*_i|.
\eeqn
Then for all $n \geqslant n_0$, we have
\beqns
E_{n+1} \leqslant  \mu \,E_{n}.
\eeqns
Consequently $x^n$ converges to $x^*$ exponentially.
\end{theorem} 

\begin{remark}
\label{notice}
{\rm Notice that if $\gamma$ is sufficiently small, e.g. 
$\gamma(1+\gamma)<\frac{2}{3}$, then for any $k<K$, there is a $\rho>0$ 
for which $\mu < 1$, so we have exponential convergence to $x^*$ whenever $x^*$ is $k$-sparse.}\hfill 
$\square$
\end{remark}
\begin{Proof}
We start with the relation \eref{wortho} with
$w=w^n$, $x^w = x^{n+1} = x^* + \eta^{n+1}$, and 
$\eta = x^{n+1}-x^* = \eta^{n+1}$, which gives
\beqns
\sum_{i=1}^N (x^*_i + \eta^{n+1}_i) \eta^{n+1}_i w^{n}_i = 0.
\eeqns
Rearranging the terms and using the fact that 
$x^*$
is supported on $T$, we get
\beqn
\label{ortho_rel}
\sum_{i=1}^N |\eta^{n+1}_i|^2 w^{n}_i 
= - \sum_{i \in T} x^*_i \eta^{n+1}_i w^{n}_i 
= - \sum_{i\in T} \frac{x^*_i}{[(x^n_i)^2 + \epsilon_n^2]^{1/2}} \eta^{n+1}_i.
\eeqn

We will prove the theorem 
by induction. Let us assume that we have
shown $E_n \leqslant   R^*$ already.
We then have, for all $i \in T$,
\beqns
|\eta^n_i| \leqslant  \| \eta^n \|_{\ell_1^N} = E_n  \leqslant  \rho |x^*_i| ~~,
\eeqns
so that
\beqn
\label{bnd1}
\frac{|x^*_i|}{[(x^n_i)^2 + \epsilon_n^2]^{1/2}}
\leqslant  \frac{|x_i^*|}{|x_i^{n}|}= \frac{|x^*_i|}{|x^*_i + \eta^n_i|} 
\leqslant  \frac{1}{1-\rho},
\eeqn
and hence \eref{ortho_rel} combined with \eref{bnd1} and NSP gives
\beqns
\sum_{i=1}^N |\eta^{n+1}_i|^2 w^{n}_i 
\leqslant \frac{1}{1-\rho} \| \eta^{n+1}_T \|_{\ell_1}
\leqslant \frac{\gamma}{1-\rho} \| \eta^{n+1}_{T^c} \|_{\ell_1}
\eeqns
At the same time, the Cauchy-Schwarz inequality combined with the above estimate yields
\begin{eqnarray}
\| \eta^{n+1}_{T^c} \|_{\ell_1}^2 
& \leqslant & \left(\sum_{i\in T^c}|\eta^{n+1}_i|^2 w^{n}_i  \right) 
\left(\sum_{i\in T^c} [(x^n_i)^2 + \epsilon_n^2]^{1/2}  \right) 
\nonumber\\
& \leqslant & \left(\sum_{i=1}^N|\eta^{n+1}_i|^2 w^{n}_i  \right) 
\left(\sum_{i\in T^c} [(\eta^n_i)^2 + \epsilon_n^2]^{1/2}  \right) 
\nonumber \\
&\leqslant& \frac{\gamma}{1-\rho} \| \eta^{n+1}_{T^c} \|_{\ell_1}
\left( \| \eta^n
\|_{\ell_1} + N \epsilon_n   \right).\label{CauSch}
\end{eqnarray}
If $\eta^{n+1}_{T^c} = 0$, then 
$x^{n+1}_{T^c}=0$. In this case $x^{n+1}$ is $k$-sparse and 
the algorithm has stopped by definition; since  
$x^{n+1}-x^*$ is in the null space $\cN$, which contains no
$k$-sparse elements other than $0$, we have already obtained
the solution $x^{n+1}=x^*$.
If $\eta^{n+1}_{T^c} \neq 0$, then after canceling the factor 
$\| \eta^{n+1}_{T^c} \|_{\ell_1}$ in \eref{CauSch}, we get
\beqns
\| \eta^{n+1}_{T^c} \|_{\ell_1}
\leqslant \frac{\gamma}{1-\rho}\left( \| \eta^n
\|_{\ell_1} + N \epsilon_n   \right)\,,
\eeqns
and thus
\beqn
\label{rec_bnd1}
\| \eta^{n+1} \|_{\ell_1} = 
\| \eta^{n+1}_T \|_{\ell_1} + 
\| \eta^{n+1}_{T^c} \|_{\ell_1} 
\leqslant 
(1+\gamma) 
\| \eta^{n+1}_{T^c} \|_{\ell_1} 
\leqslant \frac{\gamma
(1+\gamma)
}{1-\rho} 
\left ( \| \eta^n
\|_{\ell_1} + N \epsilon_n \right ).
\eeqn
Now, we also have by \eref{en} and \eref{dl11}
\beqn
\label{bnd2}
N \epsilon_n \leqslant   r(x^n)_{K+1} \leqslant \frac{1}{K+1-k}
(\|x^n - x^* \|_{\ell_1} + \sigma_k(x^*)_{\ell_1}) = 
\frac{\|\eta^n\|_{\ell_1}}{K+1-k},
\eeqn
since by assumption $\sigma_k(x^*)=0$.  This, together with \eref{rec_bnd1}, yields the desired bound,
$$E_{n+1} = \| \eta^{n+1}\|_{\ell_1} \leqslant   \frac{\gamma(1+\gamma)}{1-\rho} 
\left ( 1 + \frac{1}{K+1-k} \right ) \| \eta^{n} \|_{\ell_1}
= \mu  E_n.
$$
In particular, since
$\mu < 1$, we have $E_{n+1} \leqslant   R^*$, which completes the induction step.
It follows that  $E_{n+1} \leqslant   \mu E_n$ for all $n \geqslant  n_0$.
\end{Proof}
\begin{remark} 
{\rm Note that the precise update rule \eref{en}
for $\e_n$ does not really intervene
in this analysis. If $E_{n_0} \leqslant  R^*$,
then the estimate
\beqn
\label{contract}
E_{n+1} \leqslant   \mu_0 (E_n + N \e_n)  \, ~~~\mbox{ with } 
\mu_0 := \gamma(1+\gamma)/(1-\rho)~,
\eeqn
guarantees that all further $E_n$ will be bounded by $R^*$ as well,
provided $N\epsilon_n \leqslant (\mu_0^{-1}-1) R^*$. It is only
in guaranteeing that \eref{sball} must be satisfied for some
$n_0$ that the update rule plays a role: indeed, 
by Theorem \ref{conv}, $E_n \rightarrow 0$ for $n \rightarrow \infty$
if $\e_n$ is updated following \eref{en}, so that \eref{sball} has 
to be satisfied eventually.\\
Other update rules may work as well. If $(\e_n)_{n\in \N}$ is defined so
that it is a monotonically decreasing sequence with 
limit $\e$, then the relation \eref{contract} immediately implies that
\beqns
\limsup_{n \to \infty} E_n \leqslant  \frac{\mu_0 N \epsilon}{1-\mu_0}.
\eeqns
In particular, if $\e = 0$, then $E_n \to 0$.
The rate at which $E_n \to 0$ in this case
will depend on $\mu_0$ as well as on the rate with 
which $\e_n \to 0$. We shall not quantify this relation, except   to note 
that if $\e_n = O(\beta^n)$ for some $\beta < 1$, then 
$E_n = O(n \widetilde \mu^n)$ where $\widetilde \mu = \max (\mu_0,\beta)$.}
\hfill $\square$
\end{remark}

\subsection{Case of noisy $k$-sparse vectors}
We show here that the exponential rate of convergence to a $k$-sparse
limit vector
can be extended to the case where the ``ideal'' (i.e. $k$-sparse) target vector has been corrupted by noise and is therefore only 
``approximately $k$-sparse''.
More precisely, we no longer assume that $\cF(y)$ contains
a $k$-sparse vector; consequently the limit $\bar x$ of the $x^n$ 
need not be
an $\ell_1$-minimizer (see Theorem \ref{conv}). 
If $x$ is any $\ell_1$-minimizer in $\cF(y)$, 
Theorem \ref{conv} guarantees 
$\|\bar x -x\|_{\ell_1} \leqslant  C \sigma_k(x)_{\ell_1}$; since this is the 
best level of accuracy guaranteed in the limit, we are in this case
interested only 
in how fast $x^n$ will converge to a ball centered at $x$ with
radius given by some
(prearranged) multiple of $\sigma_k(x)_{\ell_1}$. (Note
that if $\cF(y)$ contains several $\ell_1$-minimizers, they all 
lie within a distance $C' \sigma_k(x)_{\ell_1}$ 
of each other, so that it does not matter which 
$x$ we pick.) 
We shall express the notion that $z$ is ``approximately
$k$-sparse  with gap ratio $C$'', or a 
``noisy version of a $k$-sparse vector, with gap ratio $C$'' 
by the condition 
\beqns
r(z)_k \geqslant C \sigma_k(z)_{\ell_1}
\eeqns
where $k$ is such that $\Phi$ has the NSP for some pair $K, \gamma$
such that $0\leqslant k < K-\frac{2\gamma}{1-\gamma}$ (e.g. we
could have $K=k+1$ if $\gamma < 1/2$).  If the gap ratio
$C$ is much greater than the
constant $C_1$ in \eref{est},  then 
exponential convergence can be exhibited
for a meaningful number of iterations. Note that
this class includes 
perturbations of any $k$-sparse vector 
 for which the perturbation is sufficiently small in 
$\ell^1$-norm (when compared to the unperturbed $k$-sparse vector). 

Our argument for the noisy case will closely resemble the case for 
the exact $k$-sparse vectors. However there are some 
crucial differences that justify 
our decision to separate these two cases. 

We will be interested in only the case $\e > 0$ where we recall that $\epsilon$ is the limit of the $\epsilon_n$ occurring in the algorithm, This assumption  implies 
$\sigma_k(x)_{\ell_1} > 0$, and can only happen if $x$ is not 
$K$-sparse. (As noted earlier, the
exact $k$-sparse case always corresponds to $\epsilon =0$ if
$k<K-\frac{2\gamma}{1-\gamma}$. For $k$ in the region 
$K-\frac{2\gamma}{1-\gamma}\leqslant k \leqslant K$, both $\e=0$
and $\e>0$ are theoretical possibilities.)

First, we redefine $\eta^n = x^n - x^\e$, where $x^\e$ is the minimizer
of $f_\e$ on $\cF(y)$ and $\e > 0$. 
We know from Theorem \ref{conv} that $\eta^n \to 0$. We again
set $E_n = \|\eta^n \|_{\ell_1}$.

\begin{theorem} 
\label{conv_noisy}
Given $0 < \rho < 1$, and integers $k,\,K$ with $k < K$, assume 
that $\Phi$ satisfies the NSP of order $K$
with constant 
$\gamma$ such that all the conditions of Theorem {\rm \ref{conv}} are satisfied and, in addition, 
\beqns
\mu := \frac{\gamma(1+\gamma)}{1-\rho}\left(1 + \frac{1}{K+1-k}\right) < 1.
\eeqns
Suppose $z \in \cF(y)$ 
is ``approximately $k$-sparse with gap ratio $C$'', i.e.
\beqn
\label{xcond}
r(z)_k \geqslant C \sigma_k(z)_{\ell_1} 
\eeqn
 with  $C \geqslant C_1$, where $C_1$ is as in
Theorem {\rm \ref{conv}}. 
Let $T$ stand for the set of indices of the $k$ largest entries of $x^\e$,
and $n_0$ be such that 
\beqn
\label{sball1}
E_{n_0} \leqslant  R^* := \rho \min_{i \in T} |x^\e_i| = \rho\, r(x^\e)_k.
\eeqn
Then for all $n \geqslant n_0$, we have
\beqn
\label{noisy_contract1}
E_{n+1} \leqslant  \mu E_{n} + B \sigma_k(z)_{\ell_1},
\eeqn
where $B>0$ is a constant. Similarly, if we define $\tilde E_n = 
\|x^n - z\|_{\ell_1}$, then
\beqn
\label{noisy_contract2}
\tilde E_{n+1} \leqslant  \mu \tilde E_{n} + \tilde B \sigma_k(z)_{\ell_1},
\eeqn
for $n \geqslant n_0$,
where $\tilde B > 0$ is a constant. This implies that
$x^n$ converges at an exponential (linear) rate to the ball of radius
$\widetilde B (1-\mu)^{-1} \sigma_k(z)_{\ell_1}$ centered at $z$. 
\end{theorem} 

\begin{remark}
{\rm Note that Theorem \ref{conv} trivially implies the  
inequalities \eref{noisy_contract1} and \eref{noisy_contract2} 
in the limit $n \to \infty$ since $E_n \to 0$, $\sigma_k(z)_{\ell_1} > 0$,
and $\|\bar x- z\|_{\ell_1} \leq C_1 \sigma_k(z)_{\ell_1}$. However,
Theorem \ref{conv_noisy} quantifies the event 
when it is guaranteed that the two measures of error, $E_n$ and
$\tilde E_n$, must shrink (at least) by a factor $\mu < 1$ at 
each iteration. As noted above, this corresponds to the range 
$\sigma_k(z)_{\ell_1} \lesssim  E_n, \tilde E_n \lesssim r(x^\e)_k$, 
and would be
realized if, say, $z$ is the sum of a $k$-sparse vector and a fully supported
``noise'' vector which is sufficiently small in $\ell_1$ norm. In this sense,
the theorem shows that the rate estimate of Theorem \ref{conv} extends
to a neighborhood of $k$-sparse vectors. 
}
\end{remark}

\begin{Proof}
First, note that the existence of $n_0$ is guaranteed by the fact that $E_n \to 0$
and $R^* > 0$. For the latter, note that Lemma \ref{dlemma} and Theorem
\ref{conv} imply
\beqns
r(x^\e)_k \geqslant r(z)_k - \|z-x^\e\|_{\ell_1}
\geqslant (C-C_1) \sigma_k(z)_{\ell_1},
\eeqns
so that $R^* \geqslant \rho(C-C_1)
\sigma_k(z)_{\ell_1} > 0$. 

We follow the proof of Theorem \ref{k-sparse-rate} and consider
the orthogonality relation \eref{ortho_rel}.
Since $x^\e$ is not sparse in general, we rewrite 
\eref{ortho_rel} as
\beqn
\label{ortho_rel1}
\sum_{i=1}^N |\eta^{n+1}_i|^2 w^{n}_i 
= - \sum_{i=1}^N x^\e_i \eta^{n+1}_i w^{n}_i 
= - \sum_{i\in T \cup T^c} 
\frac{x^\e_i}{[(x^n_i)^2 + \epsilon_n^2]^{1/2}} \,\eta^{n+1}_i.
\eeqn
We deal with the contribution on $T$ in the same way as before:
\beqns
\left | \sum_{i\in T} 
\frac{x^\e_i}{[(x^n_i)^2 + \epsilon_n^2]^{1/2}} \eta^{n+1}_i \right |
\leqslant \frac{1}{1-\rho} \| \eta^{n+1}_{T} \|_{\ell_1}
\leqslant \frac{\gamma}{1-\rho} \| \eta^{n+1}_{T^c} \|_{\ell_1}
\eeqns
For the contribution on $T^c$, note that 
\beqns
\beta_n := \max_{i \in T^c} \frac{|\eta^{n+1}_i|}
{[(x^n_i)^2 + \epsilon_n^2]^{1/2}}
\leqslant  \e^{-1} \| \eta^{n+1} \|_{\ell_\infty}.
\eeqns
Since $\eta^n \to 0$ we
have $\beta_n \to 0$. 
It follows that
\beqn
\label{eq625}
\left | \sum_{i\in T^c} 
\frac{x^\e_i}{[(x^n_i)^2 + \epsilon_n^2]^{1/2}}\, \eta^{n+1}_i \right |
\leqslant \beta_n \sigma_k(x^\e)_{\ell_1}
\leqslant \beta_n (\sigma_k(z)_{\ell_1}+\|x^\e-z\|_{\ell_1})
\leqslant C_2 \beta_n \sigma_k(z)_{\ell_1},
\eeqn
where the second inequality is due to Lemma \ref{dlemma}, 
the last one 
to Theorem \ref{conv}, and $C_2 = C_1+1$.
Combining these two bounds, we get
\beqns
\sum_{i=1}^N |\eta^{n+1}_i|^2 w^{n}_i 
\leqslant \frac{\gamma}{1-\rho} \| \eta^{n+1}_{T^c} \|_{\ell_1}
+ C_2 \beta_n \sigma_k(z)_{\ell_1}
\eeqns
We combine this again with a Cauchy-Schwarz estimate, to obtain
\begin{eqnarray}
\| \eta^{n+1}_{T^c} \|_{\ell_1}^2 
& \leqslant & \left(\sum_{i\in T^c} |\eta^{n+1}_i|^2 w^{n}_i\right)
\left(\sum_{i\in T^c} [(x^n_i)^2 + \e_n^2]^{1/2} \right) \nonumber \\
& \leqslant &  \left(\sum_{i=1}^N |\eta^{n+1}_i|^2 w^{n}_i \right) 
\left(\sum_{i\in T^c} [|\eta^n_i| +|x^{\e}_i| +  \e_n] \right) \nonumber \\
& \leqslant &  
\left( \frac{\gamma}{1-\rho} \| \eta^{n+1}_{T^c} \|_{\ell_1}
+ C_2 \beta_n \sigma_k(z)_{\ell_1} \right)
( \| \eta^n_{T^c} \|_{\ell_1} + \sigma_k(x^{\e})_{\ell_1} + N\e_n )\nonumber \\
& \leqslant &  
\left( \frac{\gamma}{1-\rho} \| \eta^{n+1}_{T^c} \|_{\ell_1}
+ C_2 \beta_n \sigma_k(z)_{\ell_1} \right)
( \| \eta^n_{T^c} \|_{\ell_1} + C_2 \sigma_k(z)_{\ell_1} + N\e_n )~,
\label{eq611}
\end{eqnarray}
It is easy to check that if $u^2 \leqslant Au+B$, where $A$ and $B$ are positive,
then $u \leqslant A + B/A$. Applying this to $u = \| \eta^{n+1}_{T^c} \|_{\ell_1}$
in the above estimate, we get
\beqn
\label{contract2}
\| \eta^{n+1}_{T^c} \|_{\ell_1}
\leqslant \frac{\gamma}{1-\rho} \,[ \| \eta^n_{T^c} \|_{\ell_1} + 
C_2\sigma_k(z)_{\ell_1} +N\e_n ]
+ C_3 \beta_n \sigma_k(z)_{\ell_1},
\eeqn
where $C_3 = C_2(1-\rho)/\gamma$.
Similar to \eref{bnd2}, we also have, by combining \eref{dl11} with 
(part of) the chain of inequalities \eref{eq625},
\beqn
\label{bnd3}
N \epsilon_n \leqslant  r(x^n)_{K+1} \leqslant \frac{1}{K+1-k}
(\|x^n - x^\e \|_{\ell_1} + \sigma_k(x^\e)_{\ell_1}) \le
\frac{1}{K+1-k} \left(\|\eta^n\|_{\ell_1} +
C_2 \sigma_k(z)_{\ell_1} \right),
\eeqn
and consequently \eref{contract2} becomes
\begin{eqnarray}
\| \eta^{n+1} \|_{\ell_1} & \leqslant &(1+\gamma) \| \eta^{n+1}_{T^c} \|_{\ell_1}
\label{contract3}\\
&\leqslant &\frac{\gamma(1+\gamma)}{1-\rho} \left(1+\frac{1}{K+1-k}\right)
\| \eta^n \|_{\ell_1}
+ (1+\gamma)(C_3 \beta_n + C_4) \,\sigma_k(z)_{\ell_1}\,,\nonumber 
\end{eqnarray}
where $C_4 = C_2 \gamma (1-\rho)^{-1} (1+1/(K+1-k))$.
Since the $\beta_n$ are bounded, this gives
\beqns
E_{n+1} \leqslant  \mu E_n + B \sigma_k(z)_{\ell_1}.
\eeqns
It then follows that if we pick $\widetilde \mu$ so that
$1> \widetilde \mu > \mu$, and consider the range
of $n > n_0$ such that
$E_n \geqslant (\widetilde \mu - \mu)^{-1}B 
\sigma_k(z)_{\ell_1} =: r^*$, then
\beqns
E_{n+1} \leqslant \widetilde \mu E_n.
\eeqns
Hence we are guaranteed exponential decay
of $E_n$ as long as $x^n$ is sufficiently far from its limit.
The smallest possible value of $r^*$
corresponds to the case $\widetilde \mu \approx 1$. 

To establish a rate of convergence to a comparably-sized ball 
centered at $z$, we consider $\widetilde E_n = \|x^n - z\|_{\ell_1}$. It then
follows that 
\begin{eqnarray}
\widetilde E_{n+1} 
& \leqslant  & \|x^{n+1} - x^\e\|_{\ell_1} + 
\|x^\e - z\|_{\ell_1} \cr
& \leqslant  & \mu \|x^n - x^\e\|_{\ell_1} + B \sigma_k(z)_{\ell_1} +
C_1 \sigma_k(z)_{\ell_1} \cr
& \leqslant  & \mu \|x^n - z\|_{\ell_1} + B \sigma_k(z)_{\ell_1} +
C_1(1+\mu) \sigma_k(z)_{\ell_1} \cr
& = & \mu \widetilde E_n + \tilde B \sigma_k(z)_{\ell_1},
\end{eqnarray}
which shows the claimed exponential decay and also that 
$$\limsup_{n\to \infty} \widetilde E_n \leqslant  \widetilde B (1-\mu)^{-1} 
\sigma_k(z)_{\ell_1}.$$
\end{Proof}

\section{Beyond the convex case: $\ell_\tau$-minimization for $\tau<1$}
\label{sectnonconv}

If $\Phi$ has the
NSP of order  $K$ with $\gamma<1$, then (see \S \ref{ripsect})
$\ell_1$-minimization recovers   $K$-sparse solutions to 
$\Phi x = y$ for any $y \in \R^m$ that admits such a $k$-sparse solution, 
i.e., $\ell_1$-minimization gives also $\ell_0$-minimizers, provided their
support has size at most $k$. 
In \cite{grni07}, Gribonval and Nielsen showed that in this
case,  $\ell_1$-minimization 
also gives the $\ell_\tau$-minimizers, i.e., 
$\ell_1$-minimization  also solves non-convex optimization problems of 
the type
\beqn
\label{ltau}
x^* = \argmin_{ z \in \mathcal F(y)} \| z \|_{\ell_\tau^N}^\tau, \mbox{ for } 0 <  \tau <1.
\eeqn

Let us  first  recall the results of \cite{grni07} that are of most interest 
to us here, reformulated for our setting and notations.
\begin{lemma}\emph{(\cite[Theorem 2]{grni07}).}
\label{tl1}
Assume that $x^*$ is a   $K$-sparse vector in $\mathcal F(y)$ and that $0< \tau \leqslant  1$. If
$$
\sum_{i \in T} |\eta_i|^\tau < \sum_{i \in T^c} |\eta_i|^\tau~,~
\mbox{ or, equivalently, }~
\sum_{i \in T} |\eta_i|^\tau < \frac{1}{2} \sum_{i=1}^N |\eta_i|^\tau~,
$$
for all $\eta \in \mathcal N$ and for all $T \subset \{1,\dots,N\}$ 
with $  \# T \leqslant  K$, then
$$
x^* = \argmin_{ z \in \mathcal F(y)} \| z \|_{\ell_\tau^N}^\tau.
$$
\end{lemma}
\begin{lemma}\emph{(\cite[Theorem 5]{grni07}).}\label{tl2}
Let $z \in \mathbb R^N$, $0<\tau_1 \leqslant  \tau_2 \leqslant  1$, and $K\in \mathbb N$. Then
$$
\sup_{T \subset \{1,\dots,N\}, \# T \leqslant   { K}}\frac{\sum_{i \in T} |z_i|^{\tau_1}}{\sum_{i=1}^N |z_i|^{\tau_1}} \leqslant  \sup_{T \subset \{1,\dots,N\}, \# T \leqslant  K}\frac{\sum_{i \in T} |z_i|^{\tau_2}}{\sum_{i=1}^N |z_i|^{\tau_2}}.
$$
\end{lemma}

Combining these two lemmas with the observations in  \S \ref{ripsect} leads immediately to the following result.

\begin{theorem} Fix any $0 < \tau \leqslant  1$.
If $\Phi$ satisfies the NSP of order $K$ with constant $\gamma$ then 
\begin{equation}
\label{tNSP}
\sum_{i \in T} |\eta_i|^\tau <\gamma \sum_{i \in T^c} |\eta_i|^\tau,
\end{equation}
for all $\eta \in \mathcal N$ and for all $T \subset \{1,\dots,N\}$ 
such that $\# T \leqslant  K$. 

In addition, if $\gamma <1$, and 
if there exists a $K$-sparse vector in $\mathcal F(y)$, then
this $K$-sparse vector is  the unique minimizer in $\cF(y)$ of 
$\| \cdot \|_{\ell_\tau}$.
\end{theorem}

At first sight, these results  suggest there is nothing to be gained
by carrying out $\ell_{\tau}$- rather than $\ell_1$-minimization;
in addition sparse recovery via the non-convex problems \eref{ltau}
is much harder than the more easily solvable convex relaxation problem
of $\ell_1$-minimization.

Yet, we shall show in this section that $\ell_\tau$-minimization has
unexpected benefits, and that it may be  both useful {\em and}
practically feasible via an IRLS approach. Before we start, it is 
 expedient to  introduce the following definition: 
we shall say that $\Phi$ has the 
{\it $\tau$-Null Space Property} ($\tau$-NSP) of order $K$ with constant $\gamma > 0$ if, 
for all sets $T$ of cardinality at most $K$ and all $\eta \in \cN$,
\beqn
\label{tauNSP}
\|\eta_T\|_{\ell_\tau^N}^\tau \leqslant \gamma\|\eta_{T^c}\|_{\ell_\tau^N}^\tau~.
\eeqn

In what follows we shall construct an IRLS algorithm for $\ell_\tau$- minimization.  We shall see that
\begin{itemize}
\item[(a)]  In practice, $\ell_\tau$-minimization can 
be carried out by an IRLS algorithm.
Hence, the non-convexity does not necessarily
make the problem intractable; 
\item[(b)]  In particular, if $\Phi$ satisfies the $\tau$-NSP of order $K$, and if there exists a $k$-sparse vector $x^*$ in $\mathcal F(y)$, with $k \le K-\kappa$  for suitable $\kappa$ given below, then the IRLS algorithm   converges to the 
$\ell_\tau$-minimizer $x^\tau$, which, therefore, will coincide with $x^*$; 
\item[(c)]   Surprisingly the rate of local convergence of the algorithm is superlinear; the rate is larger for 
smaller $\tau$,  increasing to approach a quadratic regime as $\tau \to 0$.  
More precisely, we will show that the local error $E_n := \| x^n - x^*\|_{\ell_\tau^N}^\tau$ satisfies
\begin{equation}
\label{ratet}
E_{n+1}  \leqslant  \mu(\gamma,\tau) E_n^{2-\tau},
\end{equation}
where $ \mu(\gamma,\tau) <1$ for $\gamma>0$ sufficiently small. 
The validity of \eref{ratet} is restricted to $x^n$ in a (small) ball centered
at $x^*$. In particular, if $x^0$ is close enough to $x^*$ then 
\eref{ratet} ensures the convergence of the algorithm { to the $k$-sparse solution $x^*$}.
\end{itemize}

{\mnew Some of these virtues of $\ell_\tau$-minimization were recently highlighted by Chartrand and his collaborators \cite{ch07,ch08,chst08}. 
Chartrand and Staneva \cite{chst08} 
 give a fine analysis of the RIP from which they can conclude   that $\ell_\tau$-minimization not only recovers $k$-sparse vectors, but that the range of $k$ for
which this recovery works is larger for smaller $\tau$.  Namely, for random Gaussian matrices, they prove that with high probability on the draw of the matrix   sparse recovery by $\ell_{\tau}$-minimization works for  
$k \leq m [c_1(\tau) + \tau c_2(\tau) \log(N/k)]^{-1}$, where  $c_1(\tau)$ is bounded and $c_2(\tau)$ decreases to zero as    $\tau \to 0$.   In particular, the dependence of the sparsity $k$ on the number $N$ of columns vanishes for $\tau \to 0$. 
These bounds give a quantitative estimate of the improvement 
provided by $\ell_\tau$-minimization vis a vis  $\ell_1$-minimization for which the range of $k$-sparsity for having exact recovery is clearly smaller (see Figure 8.4 for a numerical illustration).}
\subsection{Some useful properties of $\ell_\tau$ spaces}
\label{subsect:ltau}
We start by listing in one proposition 
some fundamental and well-known properties of 
$\ell_\tau$ spaces for $0 < \tau \leqslant  1$. For further details we refer the reader to, e.g., \cite{DeVore}.
\begin{prop}
\label{ltprop}
$~$\\
(i) Assume $0 < \tau \leqslant  1$. Then the map $z \mapsto \| z\|_{\ell_\tau^N}$ defines a quasi-norm for $\mathbb R^N$, in particular the triangle inequality holds up to a constant, i.e.,
\beqn
\label{triawc}
\| u+v\|_{\ell_\tau^N} \leqslant  C(\tau)\left( \| u\|_{\ell_\tau^N} + \| v\|_{\ell_\tau^N} \right ), \mbox{ for all } u,v \in \mathbb R^N.
\eeqn
If one considers the $\tau$-th powers of the ``$\tau$-norm'', then 
one has the { so-called} ``$\tau$-triangle inequality'':
\begin{equation}
\label{ttria}
\| u+v\|_{\ell_\tau^N}^\tau \leqslant  \| u\|_{\ell_\tau^N}^\tau  + \| v\|_{\ell_\tau^N}^\tau, \mbox{ for all } u,v \in \mathbb R^N.
\end{equation}
\\
(ii) We have, for any $0 < \tau_1 \leqslant  \tau_2 \leqslant  \infty$ 
\begin{equation}
\label{nest}
\| u\|_{\ell_{\tau_2}} \leqslant  \| u\|_{\ell_{\tau_1}}, \quad \mbox{ for all } u \in \mathbb{R}^N.
\end{equation}
We will refer to this norm estimate by writing the embedding relation 
$\ell_{\tau_1}^N  \hookrightarrow \ell_{\tau_2}^N$. \\
(iii) {\rm (}Generalized H\"older inequality{\rm)} For $0 < \tau \leqslant  1$ and $0 <p,q < \infty$ such that $\frac{1}{\tau} = \frac{1}{p}+\frac{1}{q}$, and for a positive weight  vector $w=(w_i)_{i=1}^N$ we have
\begin{equation}
\label{hoeld}
\| (u_i v_i)_{i=1}^N\|_{\ell_\tau^N(w)} \leqslant  \| u\|_{\ell_p^N(w)}\| v\|_{\ell_q^N(w)}, \mbox{ for all } u,v \in \mathbb R^N,
\end{equation}
where $\| v\|_{\ell_r^N(w)} := \left ( \sum_{i=1}^N |v_i|^r w_i \right )^{1/r}$, as usual{ , for $0<r<\infty$.}
\end{prop}

For technical reasons, it is often more convenient to employ the $\tau$-triangle inequality \eref{ttria} than \eref{triawc}; in this sense, for $\ell_\tau$-minimization $\|\cdot\|_{\ell_\tau^N}^\tau$ turns out to be more natural as a measure of error than the quasi-norm $\|\cdot\|_{\ell_\tau^N}$.\\

In order to prove the three claims (a)-(c) listed before the start of this 
subsection, we also need to generalize to
$\ell_{\tau}$ certain results previously shown only for $\ell_1$. In the following we assume $0 < \tau \leqslant  1$. We denote by 
$$
\sigma_k(z)_{\ell_\tau^N} : = \sum_{\nu > k} r(z)_\nu^\tau,
$$
the error of the best $k$-term approximation to $z$ with respect to $\|\cdot\|_{\ell_\tau^N}^\tau$. 
As a straightforward generalization of analogous results valid for the $\ell_1$-norm, we have the following two technical lemmas.
\begin{lemma}\label{tl3}
For any $j \in \{1,\dots,N\}$, we have
$$
|\sigma_j(z)_{\ell_\tau^N} -   \sigma_j(z')_{\ell_\tau^N}| \leqslant  \| z - z'\|_{\ell_\tau^N}^\tau, 
$$
for all $z,z' \in \mathbb R^N$. Moreover, for any $J>j$, we have
$$
(J-j) r(z)_j^\tau \leqslant  \sigma_j(z)_{\ell_\tau^N} \leqslant  \| z - z'\|_{\ell_\tau^N}^\tau + \sigma_j(z')_{\ell_\tau^N}.
$$
\end{lemma}

\begin{lemma}\label{tl4}
Assume that $\Phi$ has the $\tau$-NSP of order $K$   with constant $0<\gamma<1$. Then, for any $z,z' \in \mathcal F(y)$, we have 
$$
\| z' - z\|_{\ell_\tau^N}^\tau \leqslant  \frac{1+\gamma}{1-\gamma} \left (  \| z'\|_{\ell_\tau^N}^\tau -  \| z \|_{\ell_\tau^N}^\tau + 2 \sigma_K (z)_{\ell_\tau^N} \right).
$$
\end{lemma}

The proofs of these lemmas are essentially identical to the ones of 
Lemma \ref{dlemma} and Lemma \ref{balllemma}, except for substituting 
$\| \cdot \|_{\ell_\tau^N}^\tau$ for $\| \cdot \|_{\ell_1^N}$ and 
$\sigma_k (\cdot)_{\ell_\tau^N}$ for  $\sigma_k (\cdot)_{\ell_1^N}$ respectively.
\subsection{An IRLS algorithm for $\ell_\tau$-minimization}
To define an IRLS algorithm promoting $\ell_\tau$-minimization for a generic $0 < \tau \leqslant  1$, we first define a $\tau$-dependent functional
$\cJ_\tau$, generalizing $\cJ$:
\beqn
\label{defJt}
\cJ_\tau(z,w,\epsilon):= \frac{\tau}{2} \left [ \sum_{j=1}^N z_j^2w_j
+\sum_{j=1}^N\left(\epsilon^2w_j+\frac{2-\tau}{\tau} 
\frac{1}{w_j^{\frac{\tau}{2-\tau}}}\right)\right ],\quad z\in \R^N, w \in \R^N_+, \epsilon \in \R_+.
\eeqn
The desired algorithm is then defined simply by substituting $\cJ_\tau$ 
for $\cJ$ in Algorithm \ref{alg1}, keeping the same update rule \eref{en}  for $\e$. In particular we have
$$
w_j^{n+1} = \left ( (x_j^{n+1})^2 + \epsilon_{n+1}^2 \right)^{-\frac{2-\tau}{2}}, \quad j=1,\dots, N,
$$
and
$$
\cJ_\tau(x^{n+1},w^{n+1},\epsilon_{n+1}) = \sum_{j=1}^N \left ( (x_j^{n+1})^2 + \epsilon_{n+1}^2 \right)^{\frac{\tau}{2}}.
$$

Fundamental properties of the algorithm are derived in the same way as before. In particular, the values $\cJ_\tau(x^{n},w^{n},\epsilon_{n})$ decrease monotonically,
$$
\cJ_\tau(x^{n+1},w^{n+1},\epsilon_{n+1}) \leqslant  \cJ_\tau(x^{n},w^{n},\epsilon_{n}), \quad n \geqslant 0,
$$
and the iterates are bounded,
$$
\|x^n\|_{\ell_\tau^N}^\tau \leqslant  \cJ_\tau(x^{1},w^{0},\epsilon_{0}):=A_0.
$$
As in Lemma \ref{lemmawb}, the weights are uniformly bounded from below, i.e.,
$$
w_j^n \geqslant \tilde A_0,  \quad j=1,\dots,N.
$$
Moreover,  using $\cJ_\tau$ 
for $\cJ$ in Lemma \ref{boundlemma}, we can again prove the {\mnew \it asymptotic 
regularity} of the iterations, i.e.,
$$
\lim_{n \to \infty} \| x^{n+1} - x^n\|_{\ell_2^N} =0.
$$

{\mnew The first significant difference with the $\ell_1$-case arises when $\epsilon = \lim_{n \to \infty} \epsilon_n>0$. In this latter situation, we need to consider the function 
\begin{equation}
\label{targf}
f_\epsilon^\tau(z) := \sum_{j=1}^N (z_j^2 + \epsilon^2)^{\frac{\tau}{2}}.
\end{equation}
We denote by $\cZ_{\epsilon,\tau}(y)$ its set of minimizers on $\cF(y)$(since $f_{\epsilon,\tau}$ is no longer convex it may have more than one minimizer).
Even though every minimizer $z\in \cZ_{\epsilon,\tau}(y)$  still satisfies 
$$
\langle z , \eta \rangle_{w } = 0, \quad \mbox{ for all } \eta \in \cN,
$$
where $w=w^{\epsilon,\tau,z}$ is defined by  $w^{\epsilon,\tau,z}_j = ( ( z_j)^2 + \epsilon^2)^\frac{\tau-2}{2}$, $j=1,\dots,N$, the converse need no longer be true.}

{\mnew The following theorem summarizes the convergence properties on the algorithm in the case $\tau<1$.}
\begin{theorem}
\label{tconv}
{\mnew Fix $y \in \mathbb R^N$. Let $K$ (the same index as in the update rule
{\rm \eref{en} }) be chosen so that $\Phi$ satisfies the $\tau$-NSP of order $K$ with a constant $\gamma$ such that $\gamma <  1-\frac{2}{K+2}$. Let $\bar\cZ_{\epsilon,\tau}(y)$ be the set of  accumulation  points of $(x^n)_{n \in \mathbb N}$, and define $\epsilon := \lim_{n \to \infty} \epsilon_n$. Then, the algorithm has the following properties:\\
{\rm (}i{\rm )} If $\epsilon =0$, then $\bar\cZ_{\epsilon,\tau}(y)$ consists of a single point $\bar x$,  the $x^{(n)}$ converge to $\bar x$, and $\bar x$ is an $\ell_\tau$-minimizer in $\cF(y)$  which is also $K$-sparse. \\
{\rm (}ii{\rm)}
If $\epsilon >0$, then for each $\bar x\in\bar \cZ_{\epsilon,\tau}(y)$ we have
$\langle \bar x , \eta \rangle_{w^{\epsilon,\tau,\bar x}} = 0$, 
for all $\eta \in \cN$. \\ 
{\rm (}iii{\rm )}  If $z\in\cF(y)$ and $\bar x\in\bar\cZ_{\epsilon,\tau}(y)\cap \cZ_{\epsilon,\tau}(y)$,   we have
$$
\|z- \bar x  \|_{\ell_\tau^N}^\tau \leqslant  C_2 \sigma_k(z)_{\ell_\tau^N}\,,
$$
for all
$k <  K - \frac{2\gamma}{1-\gamma}$. \\
}
\end{theorem}

The proof of this theorem uses Lemmas \ref{tl1}-\ref{tl4} and follows the same arguments as for Theorem \ref{conv}.

\begin{remark}
{\rm
Unlike Theorem \ref{conv}, Theorem \ref{tconv} does not ensure that the IRLS algorithm converges to the sparsest or to the minimal $\ell_\tau$-solution. It does provide conditions that are verifiable a posteriori (e.g., $\epsilon = \lim_{n \to \infty} \epsilon_n=0$) for such convergence. The reason for this weaker result is the non-convexity of $f^\tau_\epsilon$. {\mnew (In particular, it might happen that $x^{\epsilon,\tau}$ is a local minimizer of  $f^\tau_\epsilon$, but not a global one, and the estimate in (iii) does not necessarily hold.)}  
Nevertheless, as is often the case for non-convex problems, we can establish a local convergence result that also highlights the rate we can expect for such convergence. This is the content of the following section; it will
be followed by numerical results that dovetail nicely with the theoretical results. }
\end{remark}

\subsection{Local superlinear convergence}

Throughout this section, we assume that there exists a $k$-sparse vector $x^*$ in $\cF(y)$.   We define the error vectors $\eta^n = x^n - x^* \in \cN$; we now measure the error by $\|\cdot\|_{\tau}^\tau$\,:
$$
E_n:= \| \eta^n\|_{\ell_\tau^N}^\tau.
$$ 

\begin{theorem}
\label{superth}
Assume that $\Phi$ has the $\tau$-NSP of order $ K$   with constant $\gamma \in (0,1)$ and that $\cF(y)$ contains a $k$ sparse vector $x^*$ with {\mnew $k\leq K$}. (Here $K$ is the same 
as in the definition of $\epsilon_n$ in the update 
rule {\rm \eref{en}} in Algorithm {\rm \ref{alg1}}.) Suppose that, for a given $0 < \rho < 1$, we have 
\begin{equation}
\label{tball}
E_{n_{{\!\,}_0}} \leqslant  R^*:=\left[ \rho \, r(x^*)_k \right]^\tau
\end{equation}
and define
$$
\mu :=\mu (\rho,K,\gamma,\tau,N) = 2^{1-\tau} \gamma {\mnew (1 + \gamma)} 
   A^\tau  \left ( 1 + \left ( \frac{N^{1-\tau}}{K+1-k}\right )^{2 - \tau} \right ),\quad A:= \left( r(x^*)_k^{1-\tau}
(1-\rho)^{2 - \tau} \right )^{-1}.
$$
If $\rho$ and $\gamma$ are sufficiently small so that  
\beqn
\label{assumneed}
{\mnew \mu (R^*)^{1-\tau}=} \mu\rho^{\tau(1-\tau)}r(x^*)_k^{\tau(1-\tau)}\leqslant 1,
\eeqn
 then for all $n \geqslant n_0$ we have
\begin{equation}
\label{super}
E_{n+1} \leqslant  \mu E_n^{2 - \tau}.
\end{equation}
\end{theorem}

\begin{Proof}
The proof is by induction on $n$.   We assume that
$E_n\le R^*$  and derive  \eref{super}.   As in the proof of Theorem \ref{k-sparse-rate}, we let $T$ denote
the support of $x^*$ and so $\#(T)=k$ and $r(x^*)_k$ is the smallest entry in $x^*$.
Following the proof of Theorem \ref{k-sparse-rate}, the first few lines are the same. The first difference is in the
following estimate, which holds for $i\in T$ and replaces \eref{bnd1},
\begin{eqnarray*}
\frac{|x_i^*|}{((x_i^n)^2+\epsilon_n^2)^{1-
\tau/2}}
 &\leqslant& \frac{|x_i^*|}{|x_i^*+\eta^n_i|^{2-\tau}} \leqslant  \frac{|x_i^*|}{(|x_i^*|(1-\rho))^{2-\tau}}\\
 &=&\frac{1}{|x_i^*|^{1-\tau}(1-\rho)^{2-\tau}} \leqslant A.
 \end{eqnarray*}
Starting with the orthogonality relation \eref{ortho_rel} and using the above inequality and  the embedding $\ell_{\tau_1}^N  \hookrightarrow \ell_{1}^N$, we obtain
$$
\sum_{i=1}^N |\eta_i^{n+1}|^2 w_i^n \leqslant  A \left (\sum_{i \in T} |\eta_i^{n+1}|^\tau \right)^{1/\tau}.
$$
We now apply the $\tau$-NSP to find
\begin{equation}
\label{2norm}
\| \eta^{n+1} \|^{2 \tau}_{\ell_2(w^n)} = \left ( \sum_{i=1}^N |\eta_i^{n+1}|^2 w_i^n \right )^\tau \leqslant  \gamma A^\tau \|\eta^{n+1}_{T^c} \|_{\ell_\tau^N}^\tau.
\end{equation}
At the same time, the generalized H\"older inequality (see Proposition \ref{ltprop} (iii)) for $p=2$ and $q=\frac{2\tau}{2- \tau}$, together with 
the above estimates, yields
\begin{eqnarray*}
\| \eta^{n+1}_{T^c} \|^{2\tau}_{\ell_\tau^N} &=& \| (|\eta_i^{n+1}| (w_i^n)^{-1/\tau})_{i=1}^N\|_{\ell_\tau^N(w^n;T^c)}^{2 \tau}\\
&\leqslant& \| \eta^{n+1} \|^{2 \tau}_{\ell_2^N(w^n)} \| ((w_i^n)^{-1/\tau})_{i=1}^N \|^{2 \tau}_{\ell_{2\tau/(2- \tau)}^N(w^n;T^c)}\\
&\leqslant&\gamma  A^\tau \|\eta^{n+1}_{T^c} \|_{\ell_\tau^N}^\tau \| ((w_i^n)^{-1/\tau})_{i=1}^N \|^{2 \tau}_{\ell_{2\tau/(2- \tau)}^N(w^n;T^c)}
\end{eqnarray*}
In other words,
\beqn
\label{otherwords}
\| \eta^{n+1}_{T^c} \|^{\tau}_{\ell_\tau^N}\leqslant \gamma  A^\tau   \| ((w_i^n)^{-1/\tau})_{i=1}^N \|^{2 \tau}_{\ell_{2\tau/(2- \tau)}^N(w^n;T^c)}.
\eeqn
Let us now estimate the weight term. By the $\frac{\tau}{2}$-triangle inequality \eref{ttria} we have
\begin{eqnarray*}
 \| ((w_i^n)^{-1/\tau})_{i=1}^N \|^{2 \tau}_{\ell_{2\tau/(2- \tau)}^N(w^n;T^c)}&=&\left (\sum_{i=1}^N (|\eta_i^n|^2+\epsilon_n^2)^{\frac{\tau}{2}} \right )^{2 - \tau}\\
  &\leqslant & \left (\sum_{i=1}^N (|\eta_i^n|^\tau+\epsilon_n^\tau) \right )^{2 - \tau} 
= \left (\sum_{i=1}^N |\eta_i^n|^\tau+ N \epsilon_n^\tau \right )^{2 - \tau}\\
&\leqslant & 2^{1-\tau} \left( \left( \sum_{i=1}^N |\eta_i^n|^\tau \right )^{2-\tau} + N^{2-\tau} \epsilon_n^{\tau (2 - \tau)} \right ).
\end{eqnarray*}
Now, an application of Lemma \ref{tl3} gives the following estimates
\begin{eqnarray*}
N^{2-\tau} \epsilon_n^{\tau (2 - \tau)} &=& N^{(1-\tau)(2-\tau)} (N^\tau \epsilon_n^\tau)^{2-\tau} \leqslant  N^{(1-\tau)(2-\tau)} (r(x^n)_{K+1}^\tau)^{2-\tau}\\
 &\leqslant& \left ( \frac{N^{1-\tau}}{K+1-k} \| x^n - x^*\|_{\ell_\tau^N}^\tau \right)^{2 - \tau} =\left ( \frac{N^{1-\tau}}{K+1-k} \right)^{2 - \tau} \left ( \| \eta^n \|_{\ell_\tau^N}^\tau \right)^{2 - \tau}.
\end{eqnarray*}
Using these estimates in \eref{otherwords} gives {\mnew
$$
\| \eta^{n+1}_{T^c} \|^{\tau}_{\ell_\tau^N}\leqslant 2^{1-\tau} \gamma
   A^\tau  \left ( 1 + \left ( \frac{N^{1-\tau}}{K+1-k}\right )^{2 - \tau} \right ) \left (\|\eta^{n} \|^{\tau}_{\ell_\tau^N} \right)^{2-\tau},
$$
and \eref{super} follows by a further application of the $\tau$-NSP (see \eref{rec_bnd1}).
}

\noindent
Because of the assumption \eref{assumneed}, we also have $E_{n+1}\leqslant R^*$ and so the induction can continue.
 \end{Proof}
\begin{remark}
{\rm In contrast to the $\ell_1$ case, we do not need $\mu<1$ to ensure that 
$E_n$ decreases. In fact, 
all that is needed for the error reduction is $\mu E_n^{1-\tau} < 1$ 
for some sufficiently large $n$. 
 In fact, $\mu$ could be quite large in cases where
the smallest non-zero component of the sparse vector is very small. 
 We have not observed this effect in   our examples; we expect that our analysis, although   apparently 
accurate in describing the rate of convergence (see section 8), 
is too pessimistic in estimating   the coefficient  $\mu$.}
\end{remark}
\section{Numerical results}
\label{sectnum}
In this section we present numerical experiments that illustrate
that the bounds derived in the theoretical analysis 
do manifest themselves in practice.

\subsection{Convergence rates}
\begin{figure}[htp]
  \centering 
\includegraphics[width=7.5 in, height=4 in]{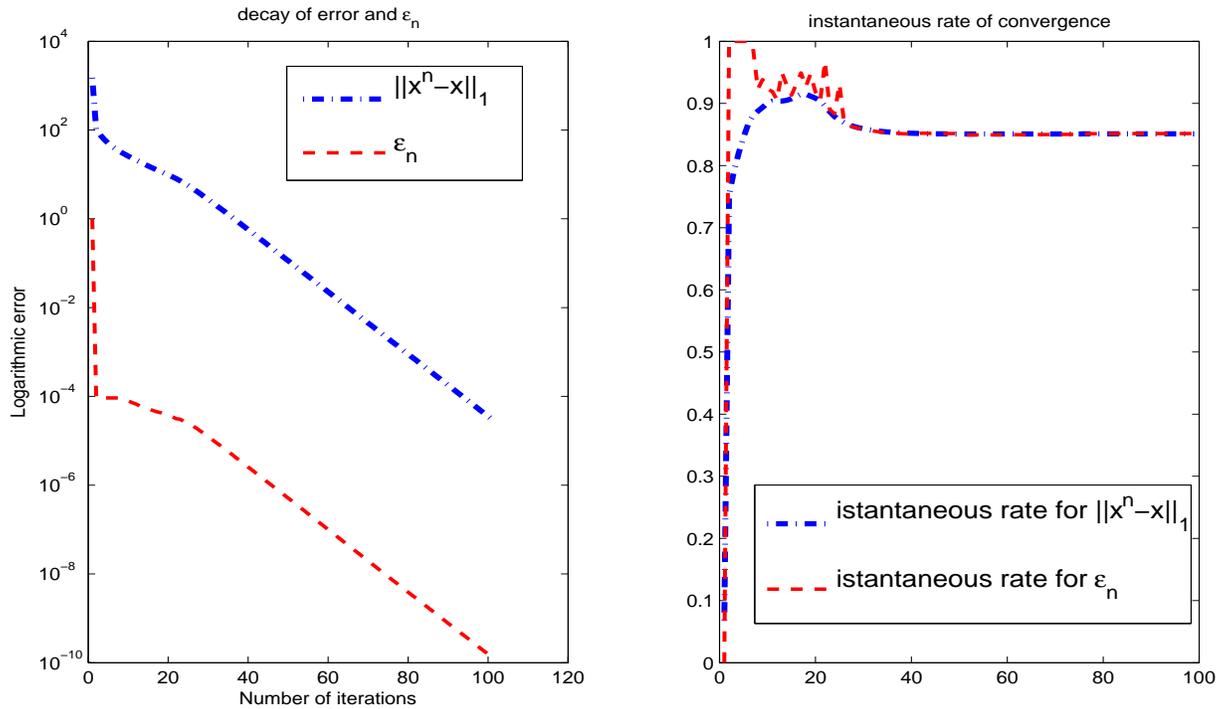}
    \caption{
{\mnew An experiment, with a matrix $\Phi$ of size $250 \times 1500$ with Gaussian
$N(0,\frac{1}{250})$ i.i.d. entries, in which recovery is sought of the 
45-sparse  vector $x^*$ represented in Figure \ref{fig2} from
its image $y=\Phi x$. Left: plot of $\log_{10}(\|x^n-x^*\|_{\ell_1})$ as a function of $n$,
where the $x^n$ are generated by Algorithm \ref{alg1}, 
with $\epsilon_n$ defined adaptively,
as in (\ref{en}). Note that the scale in the ordinate axis does not report the logarithm $0,-1,-2,\dots$, but the corresponding accuracies $10^0,10^{-1},10^{-2},\dots$ for $\|x^n-x^*\|_{\ell_1}$.
The graph also plots 
$\epsilon_n$ as a function of $n$.
Right: plot of the ratios $\|x^n-x^{n+1}\|_{\ell_1}/\|x^n-x^{n-1}\|_{\ell_1}$, 
and $(\e_n -\e_{n+1})/(\epsilon_{n-1}-\epsilon_n)$ for the same examples.}
}
\label{fig1}
\end{figure}

\begin{figure}[htp]
  \centering 
\includegraphics[width=.50\textwidth]{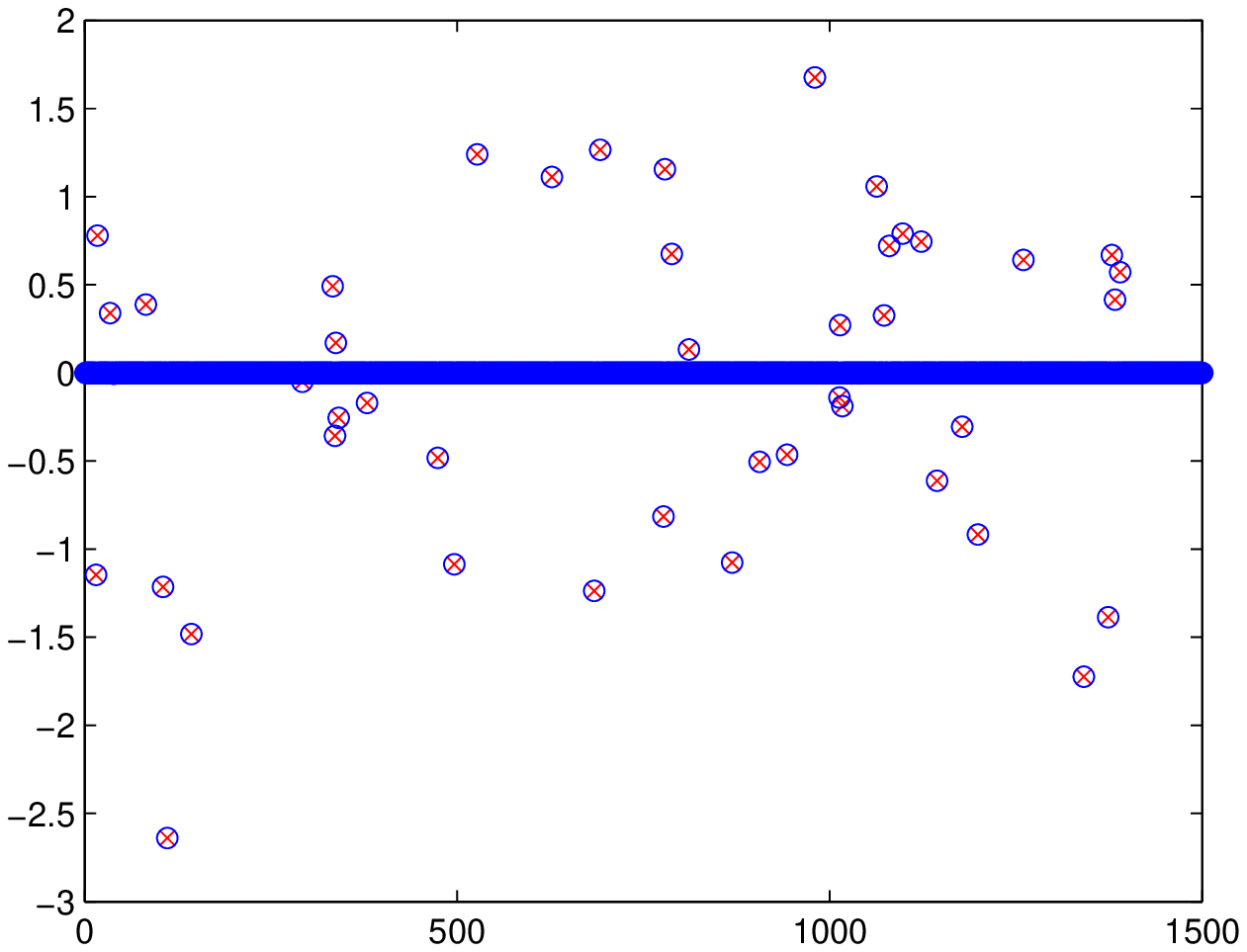}
  \caption{The sparse vector used in the example illustrated in Figure 
\ref{fig1}. This vector has dimension 1500, but only 45 non-zero entries.}
\label{fig2}
\end{figure}

We start with numerical results that confirm the linear rate of convergence of our iteratively re-weighted least square algorithm 
for $\ell_1$-minimization, and its robust recovery of sparse vectors.
In the experiments we used a matrix $\Phi$ of dimensions $m \times N$ and Gaussian 
$N(0,1/m)$  i.i.d. entries. Such matrices are known to possess 
(with high probability) the RIP property with optimal bounds \cite{BDDW,cataXX,ruveXX}.
In Figure \ref{fig1} we depict the approximation error to the unique sparsest solution shown in Figure \ref{fig2}, and the instantaneous rate of convergence. The numerical results both confirm the expected linear rate of convergence and the robust reconstruction of the sparse vector. 




Next, we compare the linear convergence achieved with $\ell_1$-minimization with the superlinear convergence obtained by the  iteratively re-weighted least square algorithm promoting $\ell_\tau$-minimization. 

\begin{figure}[ht]
  \centering 
\includegraphics[width=0.8\textwidth]{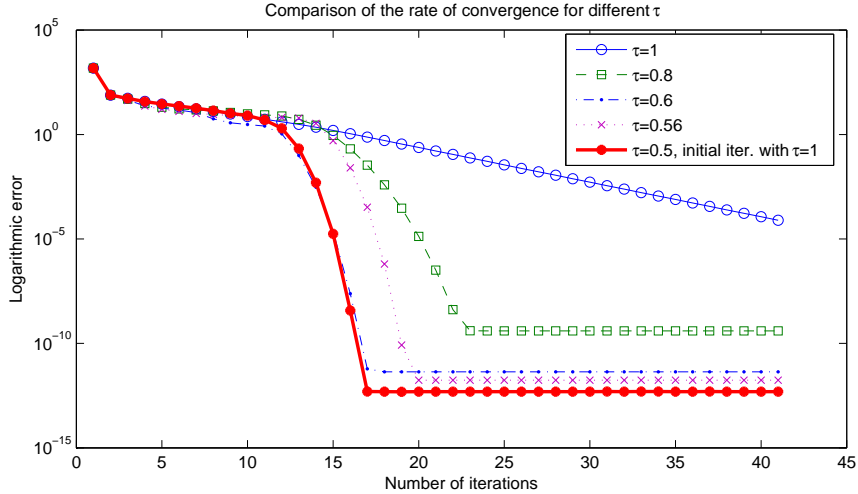}
 
  \caption{ We show the decay of logarithmic error, as a function of the
number of iterations of the algorithm for different values of $\tau$ (1,
0.8, 0.6, 0.56).
We show also the results of an experiment in which the initial $10$
iterations are performed with $\tau=1$ and the remaining iterations with
$\tau=0.5$.
}
  \label{fig3}
\end{figure}

In Figure \ref{fig3} we are interested in the comparison of the rate of convergence when our algorithm is used for different choices of $0< \tau \leqslant  1$. 
For $\tau=1, .8, .6$ and $.56$, the figure shows the error, as a function of 
the iteration step $n$, for the iterative algorithm, with different fixed values 
of $\tau$. For $\tau=1$, the rate is linear, as in
Figure \ref{fig1}. For the smaller values $\tau=.8, .6$ and $.56$ 
the iterations initially follow the same linear rate; once they
are sufficiently close to the sparse solution, the convergence rate speeds up dramatically, suggesting we have entered the region of validity of 
\eref{super}. 
For smaller values of $\tau$ numerical experiments do not always lead to
convergence: in some cases the algorithm never got to the neighborhood of the solution where convergence is ensured. However, in this case a combination of initial iterations with the $\ell_1$-inspired IRLS (for which we always have 
convergence) and later iterations with $\ell_\tau$-inspired IRLS for smaller $\tau$ 
allow again for a very fast convergence to the sparsest solution; this
is illustrated in Figure \ref{fig3} for 
the case $\tau =.5$.

\subsection{Enhanced recovery in compressed sensing and relationship with other work}
\begin{figure}[ht]
\centering 
\includegraphics[width=0.80\textwidth]{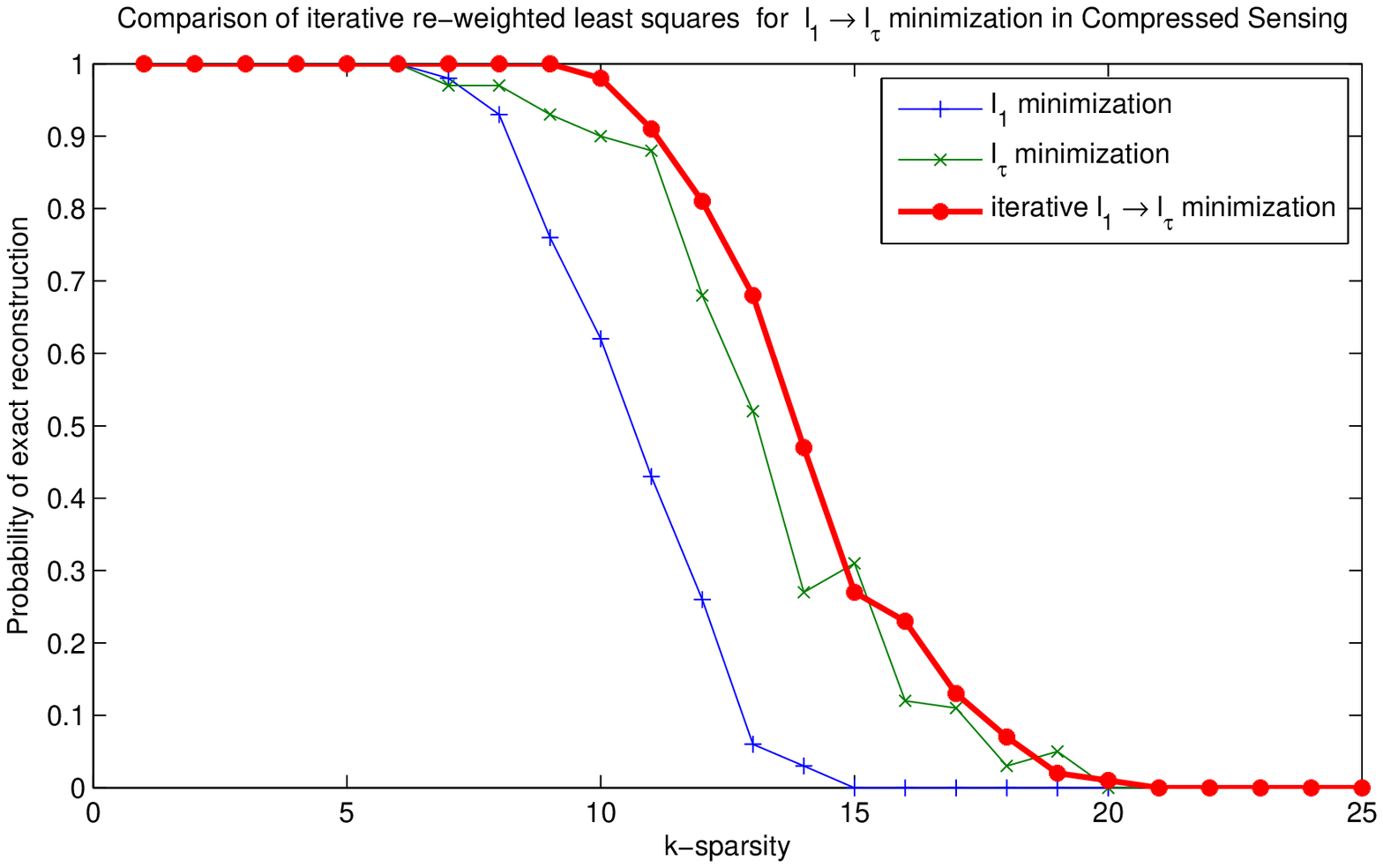}
\caption{The (experimentally determined) probability of successful
recovery of a sparse 250-dimensional vector $x$, with sparsity $k$, 
from its image $y=\Phi x$, as a function of $k$. In these experiments  
the matrix $\Phi$ is $50 \times 250$ dimensional, 
with i.i.d. Gaussian $N(0,\frac{1}{50})$ entries. 
The matrix is generated once; then, for each sparsity value $k$ shown 
in the plot, 500 attempts were 
made, for randomly generated $k$-sparse vectors $x$. 
Two different IRLS algorithms were compared, one with weights inspired by
$\ell_1$-minimization  
and the other with weights that
gradually moved from an $\ell_1$- to an $\ell_{\tau}$-minimization goal, 
with final $\tau=0.5$. We refer to \cite[Section 4]{chst08} for 
similar experiments (for different values of $\tau$), 
although realized there by fixing the 
sparsity {\it a priori} and randomly generating matrices with an 
increasing number of measurements $m$.}
  \label{fig4}
\end{figure}

   Cand\`es, Wakin, and Boyd  \cite{cawaboXX} showed, by numerical experimentation, that iteratively re-weighted $\ell_1$-minimization, with weights suggested by an $\ell_0$-minimization goal, can enhance the range 
of sparsity for which perfect reconstruction of a sparse vector ``works''
in compressed sensing. In experiments with iteratively re-weighted $\ell_2$-minimization algorithms, Chartrand and several collaborators observed 
a similar significant improvement \cite{ch07,ch08,chst08,chyi08, sachyi08}; see in particular \cite[Section 4]{chst08}; we also illustrate this in Figure \ref{fig4}.
It is to be noted that IRLS algorithms are computationally much 
less demanding than weighted $\ell_1$-minimization. 
In addition, there is, as far as we know, 
no analysis (as yet) for re-weighted $\ell_1$-minimization
that is comparable to the detailed theoretical analysis of convergence
presented here  
of our IRLS algorithm, which seems to give a realistic picture 
of the numerical computations. 
\section{Acknowledgments}
We would like to thank Yu Chen, Michael Overton  for various
conversations on the topic of this paper, and Rachel Ward for pointing out an improvement of Theorem \ref{conv}.

Ingrid  Daubechies gratefully acknowledges partial support by NSF grants DMS-0504924 and DMS-0530865.
Ronald DeVore   thanks   the Courant Institute for supporting an academic year visit  when part of 
this work was done.  He also gratefully acknowledges partial support by  Office of Naval Research Contracts
ONR-N00014-03-1-0051, ONR/DEPSCoR N00014-03-1-0675 and ONR/DEPSCoR
 N00014-00-1-0470;    the Army Research Office Contract DAAD 19-02-1-0028;  and the NSF contracts DMS-0221642 and
DMS-0200187.
Massimo Fornasier acknowledges the financial support provided by the European Union via the Individual Marie Curie fellowship MOIF-CT-2006-039438, and he thanks the Program in Applied and Computational Mathematics at Princeton University for its hospitality during the preparation of this work.
Sinan G\"unt\"urk has been supported in part by the  
National Science 
Foundation Grant CCF-0515187, an Alfred P. Sloan Research Fellowship, and an
NYU Goddard Fellowship. 
 \vskip .25in

\frenchspacing

\bibliographystyle{amsplain}
\input{DDFG14.bbl}

\end{document}

%% file: DDFG14.bbl
\providecommand{\bysame}{\leavevmode\hbox to3em{\hrulefill}\thinspace}
\providecommand{\MR}{\relax\ifhmode\unskip\space\fi MR }
\providecommand{\MRhref}[2]{%
  \href{http://www.ams.org/mathscinet-getitem?mr=#1}{#2}
}
\providecommand{\href}[2]{#2}